

\input epsf.tex

\def\2{{1\over 2}}

\def\d{\delta}
\def\a{\alpha}
\def\b{\beta}
\def\g{\gamma}

\def\s{\sigma}
\def\e{\epsilon}
\def\l{\lambda}
\def\o{\omega}

\def\fun#1#2#3{#1\colon #2\rightarrow #3}
\def\norm#1{\Vert #1 \Vert}
\def\abs#1{\vert #1 \vert}
\def\frac#1#2{{{#1} \over {#2}}}

\def\st{\;\colon\;}
\def\tends{\rightarrow}

\def\dx{\hbox{{\rm d}$x$}}

\def\dr{ {\rm d} }
\def\dt{\hbox{{\rm d}$t$}}
\def\dif{ \frac{{\rm d}}{{\rm d}t} }

\def\R{{\bf R}}
\def\N{{\bf N}}
\def\Z{{\bf Z}}

\def\T{{\bf T}}
\def\Space{{\bf S}}

\def\mod{{\rm mod}}
\def\thm#1{\vskip 1 pc\noindent{\bf Theorem #1.\quad}\sl}
\def\lem#1{\vskip 1 pc\noindent{\bf Lemma #1.\quad}\sl}
\def\prop#1{\vskip 1 pc\noindent{\bf Proposition #1.\quad}\sl}

\def\rem#1{\vskip 1 pc\noindent{\bf Remark #1.\quad}}
\def\proof{\rm\vskip 1 pc\noindent{\bf Proof.\quad}}
\def\fin{\par\hfill $\backslash\backslash\backslash$\vskip 1 pc}
\def\txt#1{\quad\hbox{#1}\quad}
\def\m{\mu}
\def\L{{\cal L}}
\def\G{\Gamma}
\def\s{\sigma}

\def\o{\omega}

\def\2{\frac{1}{2}}
\def\inn#1#2{{\langle #1 ,#2\rangle}}

\def\W{{  {\cal W}  }}
\def\V{{  {\cal V}  }}
\def\Lim{{\rm Lim}}
\def\Gcal{{  {\cal G}  }}
\def\Group{{  {\rm Group}  }}
\def\Mcal{{  {\cal M}  }}

\def\part{{\partial_{x}}}

\def\pprime{{{}^\prime{}^\prime}}



\baselineskip= 17.2pt plus 0.6pt
\font\titlefont=cmr17
\centerline{\titlefont Chaotic motions}
\vskip 1 pc
\centerline{\titlefont for a version of the Vlasov equation}
\vskip 4pc
\font\titlefont=cmr12
\centerline{         \titlefont {Ugo Bessi}\footnote*{{\rm 
Dipartimento di Matematica, Universit\`a\ Roma Tre, Largo S. 
Leonardo Murialdo, 00146 Roma, Italy.}}   }{}\footnote{}{
{{\tt email:} {\tt bessi@matrm3.mat.uniroma3.it}}} 
\vskip 0.5 pc
 
\par
\vskip 2pc
\centerline{\bf Abstract}

We consider a version of the Vlasov equation on the circle under a periodic potential $V(x,t)$ and a repulsing smooth interaction 
$W$. We suppose that the Lagrangian for the single particle has chaotic orbits; using Aubry-Mather theory and ideas of W. Gangbo, A. Tudorascu and P. Bernard, we prove that, for any initial distribution of particles, it is possible to choose their initial speed in such a way to get a chaotic orbit on $[0,+\infty)$.

\vskip 2 pc
\centerline{\bf  Introduction}
\vskip 1 pc
The Vlasov equation on the circle governs the motion of many particles on $S^1\colon=\frac{\R}{\Z}$ under the action of an external potential $V(t,x)$ and a mutual interaction $W$; we shall suppose throughout that 

\noindent $\bullet$ $V\in C^2(S^1\times S^1)$ and
$W\in C^2(S^1)$.

\noindent $\bullet$ Seen as a function on $\R$, $W$ is even: 
$W(x)=W(-x)$. Moreover, $W(x)\le 0$ and $W(x)=0$ if and only if $x\in\Z$; 
$W^\prime{}^\prime(x)<0$ when $x\in\Z$.

\vskip 1pc

Following [5], we lift the particles to $\R$; we let 
$I=[0,1)$ and we parametrize the position of the particles at time $t$ by the map $\s_t\in L^2(I,\R)$; the "Lagrangian" version of Vlasov is the ODE in $L^2(I)$
$$\left\{
\eqalign{
\ddot\s_tz&= -V^\prime(t,\s_t z)-
\int_I W^\prime(\s_tz-\s_tz^\prime)\dr z^\prime\cr
\s_0&=M\cr
\dot\s_0&=N
}
\right.     \eqno (ODE)_{Lag}$$
with $M,N\in L^2(I)$.  Note that we stick to the notation of [5] and write $\s_tz$ instead of $\s_t(z)$.

We note that $(ODE)_{Lag}$ is the Euler-Lagrange equation of the Lagrangian 
$$\fun{\L}{S^1\times L^2(I)\times L^2(I)}{\R},\qquad
\L(t,M,N)=\2\norm{N}_{L^2(I)}^2-\V(t,M)-\W(M)$$
where
$$\V(t,M)=\int_I V(t,Mz)\dr z,\qquad
\W(M)=\2\int_{I\times I}W(Mz-Mz^\prime)\dr z\dr z^\prime . $$

Now $V(t,x)$ and $W(x)$, seen as functions on $\R$, are 
$\Z$-periodic. This implies that $(ODE)_{Lag}$ has a natural invariance with respect to 
$L^2_\Z\colon=L^2(I,\Z)$: if $\s_t$ solves $(ODE)_{Lag}$ and 
$g\in L^2_\Z$, then also $\s_t+g$ solves $(ODE)_{Lag}$. Moreover, $(ODE)_{Lag}$ is also rearrangement-invariant: if 
$\s_t$ solves $(ODE)_{Lag}$ and $G$ is a measure-preserving transformation of $I$, then also $\s_t\circ G$ solves 
$(ODE)_{Lag}$. This gives us two ways to consider the Vlasov equation: either we look at it as an ODE invariant with respect to a large group of symmetries, or we concentrate on the time evolution of the {\it density} of our particles, forgetting about the labeling. Let us be a little more precise on this second approach.

Let $\fun{\pi}{\R}{S^1}$ be the natural projection; the push-forward of the Lebesgue measure on $I$ by the map 
$(\pi\circ\s_t,\dot\s_t)$ is a measure $f_t$ on 
$S^1\times\R$. If we denote by $(x,v)$ the coordinates on 
$S^1\times\R$, and by $\rho_t$ the $x$-marginal of $f_t$, then 
$f_t$ satisfies the continuity equation
$$\partial_tf_t+v\partial_xf_t=\partial_v(f_t\partial_xP_t)
\eqno (ODE)_{Meas} $$
where
$$P_t(x)=V(t,x)+\int_I W(x-x^\prime)\dr\rho_t(x^\prime) . $$
To be more precise, $f_t$ satisfies $(ODE)_{meas}$ in the weak sense, i. e.
$$-\int_{S^1\times\R}\phi\dr f_0+
\int_{[0,+\infty)}\dt\int_{S^1\times\R}
[
-\partial_t\phi-v\partial_x\phi+\partial_v\phi\partial_xP_t
]
\dr f_t  =0\qquad\forall\phi\in C^\infty_0([0,+\infty)\times S^1\times\R) . $$

In this paper, we shall follow [5] and adopt another, equivalent approach: indeed, we shall concentrate on $(ODE)_{Lag}$, but we shall quotient by the actions of $L^2_\Z$ and $\Group$.

We are interested in the relation between the Vlasov equation and the motion of a single particle, which is governed by the Lagrangian on $S^1\times S^1\times\R$
$$L(t,q,\dot q)=\2\left\vert\dot q\right\vert^2-V(t,q) . $$

Let us look at the term $-\W$ in the Lagrangian $\L$; by our hypotheses on $W$, $-\W(\s)$ is minimal if $\s z=a+bz$ with 
$a\in\R$ and $b\in L^2_\Z$; in other words, it is minimal if all the particles are grouped together on $S^1$. Thus, if 
$\fun{\s}{[0,+\infty)}{L^2(I)}$ minimizes, in some sense,
$$\int_{[0,+\infty)}\L(t,\s_t,\dot\s_t)\dt  ,  $$
we expect that the particles parametrized by $\s_t$ converge, as 
$t\tends+\infty$, to an orbit $q$ minimal for $L$. To say this precisely, we need some notation.

\vskip 1pc

\noindent {\bf Definition.} Let $K\subset\R$ be an interval; we denote by $AC(K,S^1)$ the class of absolutely continuous functions from $K$ to $S^1$. Let $c\in\R$; we say that 
$q\in AC(K,S^1)$ is $c$-minimal for $L$ if, for every 
$t_1<t_2\in K$ and 
$\tilde q\in AC((t_1,t_2),S^1)$ with $\tilde q(t_1)=q(t_1)$ and 
$\tilde q(t_2)=q(t_2)$, we have that
$$\int_{t_1}^{t_2}[L(t,q,\dot q)-c\dot q]\dt\le
\int_{t_1}^{t_2}[L(t,\tilde q,\dot{\tilde q})-c\dot{\tilde q}]\dt  .  $$
For the Lagrangian $\L$ there is a similar definition of $c$-minimality, which we postpone to the next section.

\vskip 1pc 

\noindent{\bf Definition.} Following the notation of [2], we call $\Gcal(c)$ the set of the functions 
$q\in AC(\R,S^1)$, $c$-minimal for $L$. 

\vskip 1pc

We state one of the theorems of [5]; we shall define the distance 
$dist_{weak}(M_1,M_2)$ (which is just the 2-Wasserstein distance between the measures induced by $M_1$ and $M_2$) in the next section. In the statement, with a slight abuse of notation, we identify the number $q(t)\in\R$ with the function of $L^2(I)$ constantly equal to $q(t)$.

\thm{1} Let the potentials $V$ and $W$ be as above. Let $c\in\R$ and let $\fun{M}{I}{\R}$ be monotone nondecreasing with 
$M(1-)\le M(0)+1$. Then, there is an initial speed 
$N\in L^2(I)$ such that the solution $\s_t$ of $(ODE)_{Lag}$ satisfies
$$\lim_{t\tends+\infty}\inf_{q\in\Gcal(c)}
[
dist_{weak}({\s_t},q(t))+||\dot\s_t-\dot q(t)||_{L^2(I)}
]      =0  .  $$

\rm
\vskip 1pc

Since $W^\prime(0)=0$, if $\fun{q}{\R}{S^1}$ is an orbit of $L$, then $q$ is an orbit of $\L$ too; in particular, if $L$ has chaotic orbits, so has $\L$. One could ask, however, if, for any initial distribution of particles, there are orbits, chaotic in the future with that initial distribution. In view of the precise statement, we give a few definitions.

\vskip 1pc
\noindent{\bf Definitions.}
We say that the Lagrangian $L$ admits no invariant circle of cohomology $c$ if, for some $t\in[0,1)$, the set
$$\{
q(t)\st q\in\Gcal(c)
\}   $$
is properly contained in $S^1$. We say that the interval 
$J\subset\R$ is a Birkhoff region of instability if $L$ does not admit invariant circles of cohomology $c$ for all $c\in J$.

\vskip 1pc

It is easy to see that the orbits which are in the $\a$-limit, or in the $\o$-limit of orbits in $\Gcal(c)$ are still in $\Gcal(c)$; we shall call this smaller set $\Lim(c)$.

\vskip 1pc

Our aim is to prove the following.

\thm{2} Let the interval $J\subset\R$ be a Birkhoff region of instability for $L$. Let $\fun{M}{I}{\R}$ be monotone increasing, with $M(1-)\le M(0)+1$. Let 
$\{ c_i \}_{i\in\N}\subset J$, and let 
$\{ \e_i \}_{i\in\N}\subset(0,1)$. Then there are functions 
$\fun{T_{\e_i}}{J}{\N}$ and 
$\fun{T}{J\times J}{\N}$ such that the following happens.

If 
$\{ t_i^\prime \}_{i\in\N},  \{ t_i^\pprime \}_{i\in\N}\subset(0,+\infty)$ are such that
$$t_i^\pprime-t_i^\prime\ge T_{\e_i}(c_i)
\txt{and}
t_{i+1}^\prime-t_i^\pprime\ge T(c_i,c_{i+1})$$
then there is a trajectory $\s_t$ of $(ODE)_{Lag}$ and a sequence $t_i\in(t_i^\prime,t_i^\pprime)$ such that $\s_0=M$ and
$$\inf_{q\in\Mcal(c_i)}
[
dist_{weak}(\s_{t_i},q(t_i))+
||\dot\s_{t_i}-\dot q(t_i)||_{L^2(I)}
]     \le\e_i  .  \eqno (1)$$
The sets $\Mcal(c_i)\subset\Gcal(c_i)$ will be defined in the next section. 

If, in addition, there is a class $c_\infty$ such that $c_i=c_\infty$ for large $i$, then the trajectory $\s$ has $\Lim(c_\infty)$ in the 
$\o$-limit. 

\rm
\vskip 1pc

The proof of this theorem, in section 3 below, is similar to that of theorem 2.10 (A) of [2]. In section 1, we shall recall some definitions and results from [5] and [8]; in section 2, we shall prove theorem 1. We shall do this for completeness' sake, since the proof of this theorem is distributed between [5] and [6] (see also [7]).

\vskip 2pc
\centerline{\bf \S 1}
\centerline{\bf Notation and preliminaries}

\vskip 1pc

We noticed in the introduction that $(ODE)_{Lag}$ is invariant by the actions of $L^2_\Z$ and of the group of measure-preserving maps of $I$ into itself. This prompts us to quotient $L^2(I)$ by these two groups; we recall from [5] some facts about this quotient. 

First of all, we set
$$\T\colon =\frac{L^2(I)}{L^2_\Z(I)} . $$
The space $\T$ is metric, with distance between the equivalence classes $[M]$ and $[\bar M]$ given by
$$dist_\Z([M],[\bar M])\colon=
\inf_{Z\in L^2_\Z(I)}
||
M-\bar M-Z
||_{L^2(I)}=
||
|M-\bar M|_{S^1}
||_{L^2(I)}$$
where $|x-y|_{S^1}=\min_{k\in\Z}|a-b-k|$. For any $x\in I$, we can choose measurably $Z(x)\in\Z$ such that 
$|M(x)-\bar M(x)|_{S^1}=|M(x)-\bar M(x)-Z(x)|$; this proves the second inequality above, while the first one is the definition. It also proves that the $\inf$ in the definition of $dist_\Z$ is a minimum. 

Let $\Group$ denote the group of the measurable maps of $I$ into itself which preserve Lebesgue measure and have measurable inverse; for 
$M,\bar M\in L^2(I)$ we set
$$dist_{weak}(M,\bar M)=\inf_{G\in\Group}
dist_\Z(M\circ G,\bar M) . $$
This yields that $M$ and $M\circ G$, which we would like to consider equivalent, have zero distance; however, if we say that 
$M\simeq\bar M$ when $\bar M=M\circ G$ for some 
$G\in\Group$, then the equivalence classes are not closed in 
$\T$, essentially because the $\inf$ in the definition of 
$dist_{weak}$ is not a minimum: it is possible (see [5]) that 
$dist_{weak}(M,\bar M)=0$ even if $M$ and $\bar M$ are not equivalent.  But we can consider their closure if we look at the equivalence relation from the right point of view, i. e. that of the measure induced by $M$. 

We denote by $\rm{Meas}$ the space of Borel measures on 
$S^1$, and we let $\fun{\pi}{\R}{S^1}$ be the natural projection. We introduce the map
$$\fun{\Phi}{L^2(I)}{{\rm Meas}}, \qquad
\fun{\Phi}{M}{(\pi\circ M)_\sharp\nu_0}$$
where $(\cdot)_\sharp$ denotes push-forward and $\nu_0$ is the Lebesgue measure on $I$. We note that $\Phi$ is invariant under the action of $L^2_\Z$ and $\Group$; in other words, if 
$Z\in L^2_\Z$ and $G\in\Group$, then 
$\Phi(u)=\Phi((u+Z)\circ G)$. We say that 
$M\simeq\bar M$ if $\Phi(M)=\Phi(\bar M)$. We set 
$\Space\colon=\frac{\T}{\simeq}$; on this space, we consider the metric
$$dist_\Space([M],[\bar M])=\inf\{
||M^\ast-\bar M^\ast||_{L^2(I)}   
\st
M^\ast\in[M],\quad \bar M^\ast\in[\bar M]  
\}    .   $$
The infimum above is a minimum: one can always find a minimal couple $(M^\ast,\bar M^\ast)$ with $M$ monotone and taking values in $[0,1]$, and $\bar M$ monotone and taking values in 
$[-\frac{3}{2},\frac{3}{2}]$. By lemma 2.14 of [5], $\Space$ is isometric to the space of Borel probability measures on $S^1$ with the 2-Wasserstein distance; in particular, it is a compact space. It is a consequence of proposition 2.8 of [5] that 
$dist_{weak}(M,\bar M)=dist_{{\bf S}}([M],[\bar M])$.

By proposition 2.9 of [5], which we copy below, the 
$L^2_\Z$-equivariant (or $L^2_\Z$-equivariant and 
$\Group$-equivariant) closed forms on $L^2(I)$ have a particularly simple structure.

\prop{1.1} Let $\fun{S}{L^2(I)}{\R}$ be $C^1$.

\noindent 1) If $\dr S$ is $L^2_\Z$-periodic in the sense that 
$\dr_{M+Z}S=\dr_M S$ for all $Z\in L^2_\Z(I)$, then there is a unique $C\in L^2(I)$ and a function $\fun{U}{L^2(I)}{\R}$, of class $C^1$ and $L^2_\Z$-periodic, such that
$$S(M)=U(M)+\inn{C}{M}_{L^2(I)}   .  $$

\noindent 2) If, in addition, $\fun{}{M}{\dr_M S}$ is 
rearrangement-invariant (i. e. $\dr_MS=\dr_{M\circ G}S$ for all 
$G\in\Group$), then $C$ is constant and $U$ is 
rearrangement-invariant.

\rm

\vskip 1pc

In view of the proposition above, for $c\in\R$ we define
$$\L_c(t,M,N)=\L(t,M,N)-\inn{c}{N}_{L^2(I)}  .  $$
We also define
$$L_c(t,q,\dot q)=L(t,q,\dot q)-c\dot q  .  $$
We have already defined the $c$-minimal orbits of $L$ in the introduction. In order to define the $c$-minimal orbits of $\L$, we let $K\subset\R$ be an open interval, bounded or not; following [1], we say that 
$u\in L^1_{loc}(K,L^2(I))$ is absolutely continuous (or AC) if there is $\dot u\in L^1_{loc}(K,L^2(I))$ such that, for any 
$\phi\in C^1_0(K,\R)$, we have that
$$\int_K u_t(x)\dot\phi(t)\dt=-\int_K\dot u_t(x)\phi(t)\dt . 
\eqno (1.1)$$
The equality above is in $L^2(I)$, i. e. it holds for a. e. $x\in I$; however, the exceptional set could depend on $\phi$. But it is easy to see that this is not the case, and that $\fun{}{t}{u_t(x)}$ is AC for a. e. $x\in I$.


Let $c\in\R$; we say that $\s\in AC(K,L^2(I))$ is $c$-minimal for 
$\L$ if, for any interval $[t_0,t_1]\subset K$ and any
$\tilde\s\in AC((t_0,t_1),L^2(I))$ satisfying 
$$\tilde\s_{t_1}-\s_{t_1}\in L^2_\Z(I) \txt{and} 
\tilde\s_{t_2}-\s_{t_2}\in L^2_\Z(I),$$ 
we have that
$$\int_{t_0}^{t_1}\L_c(t,\s_t,\dot\s_t)\dt\le
\int_{t_0}^{t_1}\L_c(t,\tilde\s_t,\dot{\tilde\s}_t)\dt  .  $$
Reverting to the one-particle case, we define, following [8], 
$-\a^L(c)$ as the infimum, over all the probability measures on $S^1\times S^1\times\R$ invariant by the 
Euler-Lagrange flow of $L$, of
$$\int_{S^1\times S^1\times\R}L_c(t,q,\dot q)\dr\mu(t,q,\dot q) . $$
We say that an invariant probability measure $\bar\mu$ on 
$S^1\times S^1\times\R$ is $c$-minimal if 
$$-\a^L(c)=
\int_{S^1\times S^1\times\R}L_c(t,q,\dot q)\dr\bar\mu(t,q,\dot q) . 
$$
For any $c\in\R$, there is always at least one $c$-minimal measure; we group the $c$-minimal measures in a set 
${\cal M}_{meas}(c)$. The closure of the union of all the supports of the measures in ${\cal M}_{meas}(c)$ is an invariant set; we take all the orbits which have initial condition in this set and we gather them in the set ${\cal M}(c)$; we have that
${\cal M}(c)\subset\Gcal(c)$.

Let now $n\in\N$, and let ${\cal A}_n$ be the $\s$-algebra on $I$  generated by the intervals $[\frac{i}{n},\frac{i+1}{n})$ with
$i\in(0,\dots,n-1)$; we call ${\cal C}_n$ the closed subspace of the 
${\cal A}_n$-measurable functions of $L^2(I)$, and we denote by 
$\fun{P_n}{L^2(I)}{{\cal C}_n}$ the orthogonal projection.  We have a bijection
$$\fun{D_n}{\R^n}{{\cal C}_n},\qquad
\fun{D_n}{(q_1,\dots,q_n)}{
\sum_{i=0}^{n-1}q_i1_{[\frac{i}{n},\frac{i+1}{n})}(x)
}    .   $$

\vskip 2pc

\centerline{\bf \S 2}
\centerline{\bf Proof of theorem 1}

\vskip 1pc

We shall denote by $C_{\Group}(\T)$ the class of continuous, 
$L^2_\Z$ and $\Group$-equivariant functions on $L^2(I)$; we shall also denote by 
$Mon$ the class of monotone increasing functions $M$ on 
$I=[0,1)$, such that $M(1-)\le M(0)+1$. We shall denote by 
$Mon_0$ the set of the $M\in Mon$ such that $M(0)\in[0,1]$.

\vskip 1pc

\noindent{\bf Definition of the Lax-Oleinik semigroup for} $\L$. Let 
$U\in C_{\Group}(\T)$. For $c\in\R$ and $M\in L^2(I)$, we define
$$\hat U(M)=
\inf\{
\int_0^1\L_c(t,\s_t,\dot\s_t)\dt+U(\s_1)\st 
\s\in AC([0,1],L^2(I)),\quad\s_0=M
\}   .   \eqno (2.1)  $$
Let $\a,c\in\R$; we call $\Lambda_{c,\a}$ the map
$$\fun{\Lambda_{c,\a}}{U}{\hat U+\a}. $$

\lem{2.1} The function $\hat U$ defined by (2.1) is $L^2_\Z$ and
$\Group$-equivariant.

\proof We prove that $\hat U$ is $\Group$-equivariant. Let 
$G\in\Group$; we want to show that
$$\hat U(M\circ G)=\hat U(M)  .  $$
We note that, if $\fun{}{t}{\s_t}$ is an admissible curve for the 
$\inf$ defining $\hat U(M)$, then $\s^G_t\colon =\s_t\circ G$ is admissible for the $\inf$ defining $\hat U(M\circ G)$; in other words, 
$\s^G_0=M\circ G$ and $\fun{}{t}{\s^G_t}$ is AC. Moreover, since $\L$ and $U$ are $\Group$-equivariant, we have that
$$\int_0^1\L(t,\s_t^G,\dot\s_t^G)\dt+U(\s_1^G)=
 \int_0^1\L(t,\s_t,\dot\s_t)\dt+U(\s_1)  .  $$
This implies that $\hat U(M\circ G)\le\hat U(M)$; this same formula, substituting $M\circ G$ for $M$ and $G^{-1}$ for $G$, yields the opposite inequality.

\fin

\prop{2.2} Let $U\in C_\Group(\T)$, let $M\in Mon$, and let 
$\hat U(M)$ be defined as in (2.1). Then, the following hold.

\noindent 1) $\hat U(M)$ is finite and there is 
$\s\in AC([0,1],L^2(I))$ on which $\hat U(M)$ is attained.

\noindent 2) The function $\s$ of point 1) is $c$-minimal for $\L$ on $[0,1]$ and $\s_t\in Mon$ for all $t\in[0,1]$.

\vskip 1pc
\rm

\noindent {\bf Remark.} We are going to show below that, if 
$M\in Mon$ and $\s_t$ is minimal, then $\fun{}{x}{\s_tx}$ is monotone for all $t\ge 0$; however, we only prove that $\s_t1\le\s_t0+1$ holds for a particular minimal.

\proof We need a few lemmas.

\lem{2.3} If $U\in C_\Group(\T)$, then $U$ is bounded.

\proof Let $M\in L^2(I)$ and let $\tilde M=M+Z$, where we have added $Z\in L^2_Z$ is such a way that $\tilde M$ has range in 
$[0,1]$; thus $\bar M$, the monotone rearrangement of 
$\tilde M$, belongs to $Mon_0$. Since $M$ and $\bar M$ are equivalent in the sense of section 1 (they induce the same measure on $S^1$), lemma 2.6 of [5] says that $M$ can be approximated in $L^2$ by a sequence 
$\bar M\circ G_n+Z_n$, with  
$G_n\in\Group$ and $Z_n\in L^2_{\Z}$. Since $U$ is continuous, 
$\Group$ and $L^2_\Z$-invariant, we get that $U(M)=U(\bar M)$; thus it suffices to show that $U$ is bounded on $Mon_0$. But this follows, since, in the norm topology of $L^2(I)$, $Mon_0$ is compact and $U$ is continuous. Said differently, $U$ quotients to a continuous function on ${\bf S}$; since we saw in section 1 that ${\bf S}$ is compact, we have that $U$ is bounded.

\fin

\lem{2.4} 1) For any $M\in L^2(I)$, $\hat U(M)$ is finite.

\noindent 2) If $\s$ satisfies $\s_0=M$ and
$$\int_0^1\L_c(t,\s_t,\dot\s_t)\dt+U(\s_1)\le
\hat U(M)+1  \eqno (2.2)$$
then 
$$\int_0^1||\dot\s_t||^2_{L^2(I)}\dt<C$$
for some $C>0$ independent on $U$ and $M$.

\proof We begin to prove that $\hat U(M)$ is finite. We note that 
$\V$ and $\W$ are bounded, because they are the integral of the bounded functions $V$ and $W$; $U$ is bounded by lemma 2.3. Since
$$\int_0^1[\2\norm{\dot\s_t}_{L^2(I)}^2-\inn{c}{\dot\s_t}_{L^2(I)}]\dt$$
is bounded from below, (2.1) immediately implies that 
$\hat U(M)>-\infty$.

We prove that $\hat U(M)<+\infty$. We saw in lemma 2.3 that $U$ quotients to a continuous function on the compact space ${\bf S}$. In paticular, $U$ reaches its minimum on an equivalence class 
$[\bar M]$.

Since $|a-b|_{S^1}\le 1$, the definition of $dist_\Z$ implies that the diameter of $\T$ is smaller than $1$; in particular, we can find 
$Z\in L^2_\Z$ such that, setting $\tilde M=\bar M+Z$, we have 
$||\tilde M-M||_{L^2(I)}\le 1$. We define
$$\tilde\s_t=(1-t)M+t\tilde M    $$
and we see that (2.1) implies the first inequality below; the second one follows from the fact that $\V$ and $\W$ are bounded; the last one follows from the fact that 
$\norm{\tilde M-M}_{L^2(I)}\le 1$.
$$\hat U(M)\le
\int_0^1[\2\norm{\dot{\tilde\s}_t}_{L^2(I)}^2-
\inn{c}{\dot{\tilde\s}_t}_{L^2(I)}-
\V(\tilde\s_t)-\W(\tilde\s_t)]\dt+U(\tilde\s_1)\le$$
$$\2\norm{\tilde M-M}_{L^2(I)}^2+
|c|\cdot\norm{\tilde M-M}_{L^2(I)}+C_4+U(\tilde M)\le
C_5+\min U  .  $$
This ends the proof of point 1). Let now $\s$ be as in point 2); the first inequality below is (2.2), the second one follows by the last formula.
$$\int_0^1[\2\norm{\dot\s_t}_{L^2(I)}^2-\inn{c}{\dot\s_t}_{L^2(I)}-
\V(\s_t)-\W(\s_t)]\dt\le$$
$$\hat U(M)+1-U(\s_1)\le C_5+1+\min U-U(\s_1)\le
C_5+1  .  $$
Since $\V$ and $\W$ are bounded, we get the thesis.

\fin

Following [GT], we show that we can approximate with a finite number of particles.

\lem{2.5} Let $\s$ satisfy (2.2), let $P_n$ be the projection defined in section 1 and let $\s^n_t=P_n\s_t$. Then,
$$\int_0^1\L_c(t,\s^n_t,\dot\s^n_t)\dt+U(\s^n_1)\tends
\int_0^1\L_c(t,\s_t,\dot\s_t)\dt+U(\s_1)  .  \eqno (2.3)$$

\proof Let $\s$ and $\s^n$ be as in the statement of the lemma. We assert that, as $n\tends+\infty$, 
$$\int_0^1||\s_t^n-\s_t||_{L^2(I)}\dt\tends 0
\txt{and}
\int_0^1||\dot\s^n_t-\dot\s_t||^2_{L^2(I)}\dt\tends 0 . 
\eqno (2.4)$$
The first limit in (2.2) follows  by the dominated convergence theorem: indeed, $P_n$ converges pointwise to the identity and 
$$||\s^n_t||_{L^2(I)}\le
||\s_t||_{L^2(I)}\le C_1  \qquad
\forall t\in[0,1]  .  $$
The last inequality above follows from the fact that, since $\s$ is continuous, $\s_{[0,1]}$ is compact in $L^2(I)$. The second limit in (2.3) follows analogously, since 
$\dot\s^n_t=P_n\dot\s_t\tends\dot\s_t$ and 
$\norm{\dot\s^n_t}^2_{L^2(I)}\le \norm{\dot\s_t}^2_{L^2(I)}$ for a. e. $t$; this latter function is integrable on $[0,1]$ by point 2) of lemma 2.4.

The second inequality below follows from the fact that $V$ and $W$ are Lipschitz; the last one follows from H\"older; the first formula of (2.4) implies the convergence.
$$\left\vert
\int_0^1[
\V(\s_t^n)+\W(\s^n_t)-
\V(\s_t)-\W(\s^n_t)
]\dt
\right\vert    \le$$
$$\int_0^1\dt\int_I
|
V(t,\s^n_tx)-V(\s_tx)
|\dx  +
\int_0^1\dt\int_{I\times I}
|
W(\s^n_tx-\s^n_ty)-W(\s_tx-\s_ty)
|\dx\dr y
\le$$
$$\int_0^1\dt\int_I
C_4
|
\s^n_tx-\s_tx
|\dx\le
C_4\int_0^1
||
\s_t^n-\s_t
||_{L^2(I)}\dt\tends 0  .  \eqno (2.5)$$
Applying the formula above, the second one of (2.4) and the fact that $U$ is continuous on $\T$, we get that
$$\int_0^1\L_c(t,\s^n,\dot\s^n)\dt+U(\s^n_1)=
\int_0^1[\2\norm{\dot\s_t^n}_{L^2(I)}^2-\inn{c}{\dot\s_t^n}_{L^2(I)}-
\V(\s_t^n)-\W(\s_t^n)]\dt+U(\s_1^n)\tends$$
$$\int_0^1[\2\norm{\dot\s_t}_{L^2(I)}^2-\inn{c}{\dot\s_t}_{L^2(I)}
-\V(\s_t)-\W(\s_t)]\dt+U(\s_1)=
\int_0^1\L_c(t,\s,\dot\s)\dt+U(\s_1) $$
which is the thesis.

\fin

We now introduce the finite-dimensional Lagrangian on $n$ particles, each of mass $\frac{1}{n}$:
$$L_{n,c}(t,q,\dot q)=
\frac{1}{2n}\sum_{i=1}^n|\dot q_i|^2-
\frac{c}{n}\sum_{i=1}^n\dot q_i-
\frac{1}{n}\sum_{i=1}^n V(t,q_i)-
\frac{1}{2n^2}\sum_{i,j=1}^nW(q_i-q_j)   .   $$
With this definition we have that, if the operator $D_n$ is defined as at the end of section 1,
$$L_{n,c}(t,q,\dot q)=\L_c(t,D_nq,D_n\dot q)  .  $$
This implies the following relation between the Lax-Oleinik operators of $L_n$ and $\L$:
$$\d_n\colon=\min\{
\int_0^1\L_c(t,\g_t,\dot\g_t)\dt+U(\g_1)\st
\fun{\g}{[0,1]}{{\cal C}_n},\quad\g_0=P_nM
\}  =$$
$$\min\{
\int_0^1L_{n,c}(t,q,\dot q)\dt+U(D_nq(1))\st
\fun{q}{[0,1]}{\R^n},\quad D_nq(0)=P_nM
\}  .  \eqno (2.6)$$
The minimum in the second formula above is attained by Tonelli's theorem;  clearly, if the second minimum is attained on 
$q$, the first one is attained on $D_nq$, and vice-versa.
Now, lemma 2.5 implies that 
$$\limsup_{n\tends+\infty}\d_n\le\hat U(M) . \eqno (2.7)$$
In the following, we shall show that there is $\g^n$, minimal in (2.6), is such that 
$\g^n_t\in Mon$ for $t\in[0,1]$. We shall also show that 
$\g^n$ has a subsequence $\g^{n_k}$ converging uniformly to some $\g\in AC([0,1],L^2(I))$, i. e.
$$\lim_{k\tends+\infty}\sup_{t\in[0,1]}||\g^{n_k}_t-\g_t||_{L^2(I)}=0 .\eqno (2.8)$$
We assert that this implies proposition 2.2. First of all, since 
$\g^{n_k}_t\in Mon$ and (2.8) holds, we have that $\g_t\in Mon$ for $t\in[0,1]$. Moreover, if (2.8) holds, then we see as in (2.5) that
$$\int_0^1[
\V(t,\g^{n_k}_t)+\W(\g^{n_k}_t)
] \dt\tends
\int_0^1[
\V(t,\g_t)+\W(\g_t)
] \dt   .  $$
On the other hand, the $L^2$ norm of the derivative is lower semicontinuous with respect to uniform convergence (this is true for maps valued in $\R^n$; it can be shown for maps valued in $L^2(I)$ by projecting on larger and larger subspaces), and thus
$$\int_0^1\2\norm{\dot\g_t}^2_{L^2(I)}\dt\le
\liminf_{n\tends+\infty}
\int_0^1\2\norm{\dot\g^{n_k}_t}^2_{L^2(I)}\dt  .  
$$
Moreover,
$$\int_0^1\inn{c}{\dot\s^n_t}_{L^2(I)}\dr t=
\inn{c}{\s^n_1-\s^n_0}\tends
\inn{c}{\s_1-\s_0}=
\int_0^1\inn{c}{\dot\s_t}_{L^2(I)}\dr t  .  $$
By the last three formulas and the continuity of $U$, we have that
$$\int_0^1\L_c(t,\g_t,\dot\g_t)\dt+U(\g_1)\le
\liminf_{k\tends+\infty}\left[
\int_0^1\L_c(t,\g^{n_k}_t,\dot\g^{n_k}_t)\dt+U(\g^{n_k}_1)
\right]   .  $$
This, (2.7) and the fact that $\d_{n_k}$ is attained on $\g_{n_k}$ yield
$$\int_0^1\L(t,\g_t,\dot\g_t)\dt+U(\g_1)\le\hat U(M)  .  $$
Since $\g^n_0=P_nM$, we get that $\g_0=M$; thus, by (2.1),  equality holds in the formula above; this proves part of point 1) of the thesis, namely the existence of a minimizer. The fact that 
$\hat U(M)$ is finite, is point 1) of lemma 2.4. Moreover, we get that the $\limsup$ of (2.7) is actually a limit and is equal to 
$\hat U(M)$. 

As for point 2), we have just seen that $\g_t\in Mon$ for 
$t\in[0,1]$. The fact that $\s_t$ is $c$-minimal follows immediately from the fact that $\s_t$ minimizes in (2.1).

We need the next three lemmas to prove that, when $M\in Mon$, (2.8) holds. Naturally, there are similar lemmas in [5] and [6]: "self contained" is often synonymous with "reinventing the wheel".

\lem{2.6} Let $P_nM$ be monotone increasing, and let the first 
$\min$ of (2.6) be attained on $\g^n$. Then $\fun{}{x}{\g^n_tx}$ is monotone increasing for all $t\in[0,1]$.

\proof We have seen after (2.6) that there is a function 
$(q_1(t),\dots,q_n(t))$, minimal in the second formula of (2.6),  such that 
$\g^n_t=D_n(q_1(t),\dots,q_n(t))$; since $P_nM$ is monotone increasing, we have that 
$q_1(0)\le q_2(0)\le\dots\le q_n(0)$. We must prove that, up to rearranging the indices for which $q_i(0)=q_{i+1}(0)=\dots=q_{i+l}(0)$, we have
$q_1(t)\le q_2(t)\le\dots\le q_n(t)$ for $t\in[0,1]$. Since this follows from the fact that the orbits $q_i$ are an ordered set, it suffices to prove that there are no two times 
$t_1<t_2\in[0,1]$ and two indices $i<j$ such that
$$q_i(t_1)<q_j(t_1)\txt{and}
q_i(t_2)>q_j(t_2)  .  \eqno (2.9)$$
We shall argue by contradiction, supposing that the formula above holds. Let us define a new orbit $\tilde q$ as 
$$\tilde q_i(t)=\min(q_i(t),q_j(t)),\quad
\tilde q_j(t)=\max(q_i(t),q_j(t)),\quad
\tilde q_l(t)=q_l(t)\txt{for}l\not=i,j  .  $$
Since $P_nM$ is monotone, we have that $\tilde q(0)=q(0)$; using the fact that $U$ is $\Group$-invariant, we get that
$$\int_0^1L_{n,c}(t,q,\dot q)\dt+U(D_nq(1))=
\int_0^1L_{n,c}(t,\tilde q,\dot{\tilde q})\dt+U(D_n\tilde q(1))  .
\eqno (2.10)$$
On the other hand, $\tilde q$ is not minimal. To show this, let us suppose by contradiction that $\tilde q$ is minimal; in particular, it is a solution of the Euler-Lagrange equation and it is $C^1$. Since $q$ is continuous, by (2.9) there is $\bar t\in(t_1,t_2)$ such that
$q_i(\bar t)=q_j(\bar t)=\tilde q_i(\bar t)=\tilde q_j(\bar t)$. We have seen that $\tilde q_i$ and $\tilde q_j$ are $C^1$; this and their definition imply that
$$\dot q_i(\bar t)=\dot q_j(\bar t)=\dot{\tilde q}_i(\bar t)=
\dot{\tilde q}_j(\bar t)  .  $$
Now the solutions $q$ and $\tilde q$ have the same initial conditions at $\bar t$ but they do not coincide, a contradiction. 

Since $\tilde q$ is not minimal, we can find $\hat q$ with the same boundary conditions of $\tilde q$ such that
$$\int_0^1L_{n,c}(t,\hat q,\dot{\hat q})\dt+U(D_n\hat q(1))<
\int_0^1L_{n,c}(t,\tilde q,\dot{\tilde q})\dt+U(D_n\tilde q(1))  .  $$
We recall that $\hat q(0)=\tilde q(0)=q(0)$, while 
$\hat q(1)=\tilde q(1)$ coincides with $q(1)$ up to rearranging indices, and thus $U(D_n q(1))=U(D_n\hat q(1))$; this, together with the last formula and (2.10), contradict the minimality of $q$.

\fin

\lem{2.7} Let $P_nM\in Mon$. Then there is $\g$, minimal in the first formula of (2.6), such that $\g_t\in Mon$ for all $t\in[0,1]$.

\proof We have seen in the last lemma that, if $\g_0$ is monotone, then $\g_t$ is monotone for $t\in[0,1]$; thus, we have only to prove that, if $\g_0(1)-\g_0(0)\le 1$, then there is a minimal $\g$ such that
$\g_t(1)-\g_t(0)\le 1$ for all $t\in[0,1]$.

We have also seen that there is $\fun{q}{[0,1]}{\R^n}$, minimal in the second formula of (2.6), such that $\g_t=D_nq(t)$.

There are two cases; the first one is 
$\g_0(1)-\g_0(0)<1$ or, equivalently, $q_n(0)-q_1(0)<1$.
Let us suppose by contradiction that the thesis is false; then there is $t_1\in(0,1]$ such that $\g_{t_1}(1)-\g_{t_1}(0)> 1$; since 
$\g_t=D_nq(t)$, this means that 
$$q_n(t_1)-q_1(t_1)>1  .  \eqno (2.11)$$
Let us consider the orbit
$$\tilde q_1(t)=q_1(t)+1, \quad
\tilde q_i(t)=q_i(t)\txt{for} i\not=1  .  $$
Since $L_{n,c}$ is invariant by integer translations and permutation of indices, $\tilde q$ is a minimizer in the second formula of (2.6) with initial condition $\tilde q(0)$. Since $q_n(0)-q_1(0)< 1$, we get that 
$\tilde q_1(0)>\tilde q_n(0)$, while by (2.11) we have that 
$\tilde q_1(t_1)<\tilde q_n(t_1)$. In other words, $\tilde q_1$ and 
$\tilde q_n$ cross; but this is forbidden by the last lemma. 

In the second case, $q_n(0)=q_1(0)+1$; we set 
$\tilde q_n(0)=q_n(0)-\e$ and $\tilde q_i(0)=q_i(0)$ for $i<n$. We assert that $\hat U(D_n\tilde q(0))\tends\hat U(D_n q(0))$ for
$\e\tends 0$; indeed, here we are dealing with the Lax-Oleinik semigroup in finite dimension, and it is a standard theorem that, in this case, $\hat U$ it is Lipschitz.

Let $\tilde q$ minimize in (2.6) with initial condition $\tilde q(0)$; since $\tilde q_n(0)<\tilde q_1(0)+1$ by definition, the last paragraph implies that
$D_n\tilde q(t)\in Mon$ for all $t$. We call 
$\bar q$ the limit of the orbits $\tilde q$ as $\e\tends 0$; again because we are dealing with a Lagrangian in $\R^n$, it is standard to see that this limit exists along a subsequence. It is a standard fact (which we reproved after (2.6)) that the limit $D_n\bar q$ minimizes; moreover, it belongs to $Mon$, being the limit of functions in $Mon$.

\fin

We now end the proof of (2.8); we have seen that proposition 2 follows from this formula. Let $\g^n$ be minimal in (2.6); by point 2) of lemma 2.4 and the fact that $U$ is continuous,
$$\int_0^1||\dot\g^n_t||^2_{L^2(I)}\dt\le C$$
for some $C>0$ independent on $n$. As a consequence, the 
$\2$-H\"older norm of $\fun{\g^n}{[0,1]}{L^2(I)}$ is bounded uniformly in $n$; since the $L^2$-norm of 
$\g^n_0=P_nM$ is bounded too, we get that
$$||\g^n_t||_{L^2(I)}\le C_3\quad\forall t\in[0,1]$$
for some $C_3>0$ independent on $n$. By lemma 2.7, we can suppose that $\g^n_t\in Mon$ for all $t\in[0,1]$; in particular, 
$\g^n_t(x)-\g^n_t(0)\le 1$ for all $t\in[0,1]$ and $x\in I$. Together with the last formula, this implies that
$$||\g^n_t||_{L^\infty(I)}\le C_4\quad\forall t\in[0,1] . $$
Since $\g_t\in Mon$, by the last formula we have that, for all 
$t\in[0,1]$, 
$$\g^n_t\in A\colon=Mon\cap
\{
u\in L^2(I)\st ||u||_{L^\infty(I)}\le C_4
\}  .  $$
Thus $\fun{\g^n}{[0,1]}{L^2(I)}$ is a uniformly $\2$-H\"older sequence with values in the set $A$, which is compact in $L^2(I)$ by one of Helly's theorems; now (2.8) follows by Ascoli-Arzel\`a.

\fin

We now begin to analyze the Lax-Oleinik operator 
$\Lambda_{c,\a}$, with a view to prove that, for a  suitable choice of $\a$, it has a fixed point in $C_{\Group}(\T)$.

\lem{2.8} There is $L>0$ such that the following holds. Let 
$U\in C_{\Group}(\T)$, and let $\hat U$ be defined as in (2.1). Then $\hat U$ is $L$-Lipschitz for $dist_\Z$.

\proof Let $\hat U$ be defined as in (2.1), and let 
$M_1,M_2\in L^2(I)$; we have to prove that
$$|
\hat U(M_1)-\hat U(M_2)
|\le 
L dist_\Z(M_1,M_2)  .  \eqno (2.12)$$
By section 1 and lemma 2.1, up to adding elements of $L^2_\Z$ to $M_1$ and 
$M_2$, we can suppose that $M_1$ takes values in $[0,1]$, 
$M_2$ takes values in $[-1,2]$ and 
$$dist_\Z(M_1,M_2)=||M_1-M_2||_{L^2(I)}  .  \eqno (2.13)$$
Let $\e>0$ be fixed; by (2.1), we can find $\s^1_t$ with 
$\s^1_0=M_1$ such that
$$\hat U(M_1)\ge
\int_0^1\L_c(t,\s^1_t,\dot\s^1_t)\dt+U(\s^1_1)-\e . \eqno (2.14) $$
We recall that we cannot use proposition 2.2, because $M_1$ belongs to $L^2(I)$, not to $Mon$. We set 
$$\s^2_t=\s^1_t+(1-t)(M_2-M_1) . $$
By (2.14), the definition of $\hat U(M_2)$ and the fact that 
$\s^1_1=\s^2_1$, we get the first inequality below.
$$\hat U(M_2)-\hat U(M_1)\le
\int_0^1\L_c(t,\s^2_t,\dot\s^2_t)\dt-
\int_0^1\L_c(t,\s^1_t,\dot\s^1_t)\dt+\e\le$$
$$-\int_0^1\dt\int_I\dot\s_t(x)[M_2(x)-M_1(x)]\dx+
\2\int_0^1\dt\int_I|M_2(x)-M_1(x)|^2\dx+$$
$$c\int_0^1\dt\int_I(M_2(x)-M_1(x))\dx+
2A\int_0^1(1-t)\dt\int_I  |M_2(x)-M_1(x)|\dx+\e . $$
The second inequality above comes from the fact that the potentials $V$ and $W$ are $A$-Lipschitz for some $A>0$. Now we apply Cauchy-Schwarz to the formula above; we recall that, by (2.14) and lemma 2.4, we have $||\dot\s||_{L^2([0,1]\times I)}\le C$ with $C$ independent on $U$ and $M$. We get
$$\hat U(M_2)-\hat U(M_1)\le
\tilde L(
\norm{M_1-M_2}_{L^2(I)}^2+\norm{M_1-M_2}_{L^2(I)}
)+\e  .  $$
Since $M_1$, $M_2$ take values in $[-1,2]$, the last formula implies that
$$\hat U(M_2)-\hat U(M_1)\le L||M_1-M_2||_{L^2(I)}+\e . $$
Using (2.13) and the fact that $\e>0$ is arbitrary, we get
$$\hat U(M_2)-\hat U(M_1)\le
L dist_\Z(M_1,M_2)  .  $$
Interchanging the r\^ oles of $M_1$ and $M_2$, we get (2.12).

\fin

\lem{2.9} Let $\hat U$ and $L>0$ be as in the last lemma. Then 
$\hat U$ quotients to a function on ${\bf S}$, which is 
$L$-Lipschitz for $dist_{{\bf S}}$.

\proof Let $M_1,M_2\in L^2(I)$; we know by section 1 that there are $\tilde M_1\in [M_1]$ and $\tilde M_2\in [M_2]$ such that 
$\tilde M_1$ is monotone and has range in $[0,1]$, $\tilde M_2$ is monotone and has range in $[-\2,\frac{3}{2}]$ and
$$dist_{{\bf S}}([M_1],[M_2])=
||\tilde M_1-\tilde M_2||_{L^2(I)}=
dist_\Z(\tilde M_1,\tilde M_2) . $$
By lemma 2.6 of [5], $M_1$ can be approximated in $L^2(I)$ by a sequence $(\tilde M_1+Z_i)\circ G_i$ with $Z_i\in L^2_\Z$ and $G_i\in\Group$, and the same holds for $M_2$. Since 
$dist_\Z(M,N)\le||M-N||_{L^2(I)}$, lemma 2.8 implies that $\hat U$ is continuous on $L^2(I)$; since $\hat U$ is equivariant too (lemma 2.1), we get the first equality below.
$$|
\hat U(M_1)-\hat U(M_2)
|  =
|
\hat U(\tilde M_1)-\hat U(\tilde M_2)
|\le$$
$$L dist_\Z(
\tilde M_1,\tilde M_2
)=
L dist_{{\bf S}}([M_1],[M_2]) . $$
The inequality above comes from lemma 2.8; the last formula implies the thesis. 

\fin

At the beginning of this section, we defined the operator 
$\Lambda_{c,\a}$; by the last lemma, it brings $C_{\Group}(\T)$ into itself. As usual for Lax-Oleinik operators, it is non-expansive.

\lem{2.10} The operator $\Lambda_{c,\a}$ from $C_{\Group}(\T)$ into itself is $1$-Lipschitz for the $\sup$ norm on 
$C_{\Group}(\T)$.

\proof From (2.1) we get that, if $U_1\le U_2$, then
$\hat U_1\le\hat U_2$; moreover, if $\b\in\R$, then 
$(U+\b)\hat{}=\hat U+\b$. Thus, if 
$U_3,U_4\in C_{\Group}(\T)$, we have
$$\hat U_4\le (U_3+||U_4-U_3||_{\sup})\hat{}=
\hat U_3+ ||U_4-U_3||_{\sup}  .  $$
Interchanging the r\^ oles of $U_3$ and $U_4$, we get the thesis.

\fin

\prop{2.11} There is a unique constant $\a(c)$ such that the operator $\Lambda_{c,\a(c)}$ has a fixed point in 
$C_{\Group}(\T)$. If $L$ is as in lemma 2.9, the fixed points are 
$L$-Lipschitz on ${\bf S}$ for $dist_{{\bf S}}$.

\proof  Let us denote by $Lip_L({\bf S})$ the functions $U$ defined on ${\bf S}$ which are $L$-Lipschitz for $dist_{\bf S}$. We see as in lemma 2.9 that, if $U\in C_{\Group}(\T)$ is $L$-Lipschitz for 
$\norm{\cdot}_{L^2(I)}$, then $U$ quotients to a function on 
${\bf S}$ which is $L$-Lipschitz for $dist_{{\bf S}}$; vice-versa, any function in $Lip_L({\bf S})$ lifts to an $L$-Lipschitz function in 
$C_{\Group}(\T)$: just compose with the projection on $\Space$. As a consequence, it suffices to find a fixed point of 
$\Lambda_{c,\a}$ on $Lip_L({\bf S})$.

We have seen in lemma 2.10 that, if $a\in\R$,
$$\Lambda_{c,\a}(U+a)=\Lambda_{c,\a}U+a .  $$
This prompts us to follow [4] and to define the space
$$Lip_L^0({\bf S})=\{
U\in Lip_L({\bf S})\st U([f])=0
\}  $$
where $\fun{f}{I}{S^1}$ is the function $f(x)\equiv 0$ and $[\cdot]$ denotes the equivalence class. This is a set of equicontinuous functions on ${\bf S}$; since ${\bf S}$ is bounded, 
$Lip_L^0({\bf S})$ is equibounded too. 

We set
$$\tilde\Lambda_c(U)=\hat U-\hat U([f])  .  $$
By section 1, ${\bf S}$ is  compact set for $dist_{{\bf S}}$; thus, 
$Lip^0_L({\bf S})$ is, by Ascoli-Arzel\`a, a compact set of 
$C({\bf S},\R)$; moreover, it is convex. The function 
$\tilde\Lambda_c$ brings this set into itself by lemma 2.9 and is continuous (actually, $2$-Lipschitz) by lemma 2.10; thus, 
$\tilde\Lambda_c$ has a fixed point $U_{fix}$ by the Schauder fixed point theorem. If we set 
$$\a(c)=-\hat U_{fix}([f])$$
we see by the definition of $\Lambda_{c,\a(c)}$ that
$$\Lambda_{c,\a(c)}(U_{fix})=\hat U_{fix}+\a(c)=
\hat U_{fix}-\hat U_{fix}([f])=
\tilde\Lambda_c(U_{fix})=U_{fix}$$
i. e. that $U_{fix}$ is a fixed point of $\Lambda_{c,\a(c)}$ on 
$Lip_L({\bf S})$. 

It remains to prove that $\a(c)$ is unique. Let us suppose by contradiction that there is $\a\not=\a(c)$ such that 
$\Lambda_{c,\a}$ has a fixed point $U_1$; to fix ideas, let 
$\a>\a(c)$. Clearly, $U_1$ is a fixed point of $\Lambda_{c,\a}^n$, while $U_{fix}$ is a fixed point of $\Lambda_{c,\a(c)}^n$; this implies the first equality below.
$$\norm{U_1-U_{fix}}_{\sup}=
||
\Lambda^n_{c,\a}U_1-\Lambda^n_{c,\a(c)}U_{fix}
||_{\sup}=$$
$$
||
\Lambda^n_{c,\a}U_1-\Lambda^n_{c,\a}U_{fix}+n(\a-\a(c))
||_{\sup}\ge
-||
U_1-U_{fix}
||_{\sup} +n(\a-\a(c))      $$
where the last inequality comes from lemma 2.10.
But $U_1$ and $U_{fix}$ are bounded by lemma 2.3; we let 
$n\tends+\infty$ in the formula above and, recalling that 
$\a>\a(c)$, we get a contradiction. The case $\a<\a(c)$ is analogous.

\fin

We now want to prove that the function $\a(c)$ of the last lemma coincides with the function $\a^L(c)$ we defined in section 1. The connection is provided by the following fact (see [4] for a proof): $\a^L(c)$ is the unique $\a\in\R$ for which the Lax-Oleinik operator 
$$\fun{\hat\Lambda}{C(S^1)}{C(S^1)}$$
$$(\hat\Lambda u)(x)=\a+\min\{
\int_0^1L_c(t,q,\dot q)\dt+u(q(1))\st q(0)=x
\}$$
has a fixed point. The minimum above is attained by Tonelli's theorem.

\lem{2.12} $\a^L(c)=\a(c)$.

\proof By the last proposition, there is a fixed point $U$ of 
$\Lambda_{c,\a(c)}$. Since $U$ is a fixed point of 
$\Lambda_{c,\a(c)}^k$ too, we can apply proposition 2.2 to 
$\Lambda_{c,\a(c)}^k$: for a constant function $M\equiv x_0$, we can find 
$\bar\s\in AC([0,k],L^2(I))$, $c$-minimal for $\L$ and such that
$\bar\s_0=M$ and 
$$U(M)=\int_0^k[\L_c(t,\bar\s_t,\dot{\bar \s_t})+\a(c)]\dt+
U(\bar\s_k)  .   $$
It can be proven as in [4] that $\bar\s$ does not depend on $k$, but we shall not need this in the proof.
Since $||U||_{Sup}\le C$ by lemma 2.3, we get that
$$\left\vert
\int_0^k[\L_c(t,\bar\s_t,\dot{\bar \s_t})+\a(c)]\dt
\right\vert
\le 2C \quad\forall k\in\N . \eqno (2.15)$$
On the other hand, let $u$ be a fixed point of $\hat\Lambda$ and let $\bar q$ realize the $\min$ in the definition of 
$(\hat\Lambda^k u)(x_0)$; we get as above that
$$u(x_0)=\int_0^k[L(t,\bar q,\dot{\bar q})+\a^L(c)]\dt+
u(\bar q(k))  .  $$
Since $u$ is bounded (it is a continuous function on $S^1$) we get again that, possibly increasing $C$,
$$\left\vert
\int_0^k[L_c(t,\bar q,\dot{\bar q})+\a^L(c)]\dt
\right\vert
\le 2C \quad\forall k\in\N . \eqno (2.16)$$
Since $\bar q(0)=M\equiv x_0$, (2.1) implies the first inequality below; as we did above, we identify the number $\bar q(t)$ with the function in $L^2(I)$ constantly equal to $\bar q(t)$; we could have underlined this identification using the operator $D_1$.
$$U(M)=(\Lambda^k_{c,\a(c)}U)(M)\le\int_0^k
[\L_c(t,\bar q,\dot{\bar q})+\a(c)]\dt+
U(\bar q(k))=$$
$$\int_0^k[L_c(t,\bar q,\dot{\bar q})+\a(c)]\dt+
U(\bar q(k))=
\int_0^k[L_c(t,\bar q,\dot{\bar q})+\a^L(c)]\dt+
k[\a(c)-\a^L(c)]+U(\bar q(k))\le$$
$$k[\a(c)-\a^L(c)]+U(\bar q(k))+2C  .  $$
The second equality above comes from the fact that $W(0)=0$, the second inequality from (2.16). Since $U$ is bounded, taking $k$ large in the formula above, we get that 
$$\a(c)\ge\a^L(c)  .  \eqno (2.17)$$
We prove the opposite inequality. Let $\bar\s$ be as in (2.15); by lemma 2.5, we can find $n=n(k,\bar\s)$ so large that, defining $q=(q_1,\dots,q_n)$ by $D_nq(t)=P_n\bar\s_t$, we have
$$\int_0^k\L_c(t,\bar\s_t,\dot{\bar \s_t})\dt
+U(\s_k)\ge
\int_0^kL_{n,c}(t, q,\dot{q})\dt+
U(D_n q(k))-1  .  $$
We recall that $||U||_{\sup}\le C$; we set $C_1=2C+1$ and we get from the formula above that
$$\int_0^k\L_c(t,\bar\s_t,\dot{\bar\s}_t)\dt+C_1\ge
\int_0^kL_{n,c}(t,q,\dot{q})\dt  . $$
Note that $C_1$ does not depend on $k$ and $n$. Since
$$L_{n,c}(t,q,\dot q)=\frac{1}{n}\sum_{i=1}^nL_c(t,q_i,\dot q_i)-
\frac{1}{2n^2}\sum_{i,j=1}^nW(q_i-q_j)$$
and $W\le 0$, we get that there is $i\in(1,\dots,n)$ such that
$$\int_0^k\L_c(t,\bar\s_t,\dot{\bar\s}_t)\dt+C_1\ge
\int_0^kL_c(t,q_i,\dot{q_i})\dt  .  $$
Taking $\tilde q$  minimal in the definition of 
$(\hat\Lambda^k u)(q_i(0))$, we get the inequality below.
$$u(q_i(0))=
\int_0^k[
L_c(t,\tilde q,\dot{\tilde q})+\a^L(c)
]\dt  +u(\tilde q(k))\le
\int_0^k[
L_c(t,q_i,\dot{q}_i)+\a^L(c)
]\dt  +u(q_i(k))  .   $$
Since $u$ is bounded, by the last two formulas we get
$$\int_0^k L_c(t,\tilde q,\dot{\tilde q})\dt\le
\int_0^k\L_c(t,\bar\s_t,\dot{\bar \s}_t)\dt+C_2    $$
for some $C_2>0$ independent on $n$ and $k$.
We have put into $C_2$ the bounded contribution of $u$. 
We note that (2.16) holds for $\tilde q$ too; indeed, in the proof we only used that $\tilde q$ is minimal and that $u$ is bounded. Now, (2.16) for $\tilde q$ and (2.15) yield the first inequality below; the last formula yields the second one. 
$$k[\a(c)-\a^L(c)]\le  
\int_0^k L_c(t,\tilde q,\dot{\tilde q})\dt
-\int_0^k\L_c(t,\bar\s_t,\dot{\bar \s}_t)\dt +4C\le C_3 .  $$
Letting $k\tends+\infty$, we get the inequality opposite to (2.17), and we are done.

\fin

The next lemma tells us the the sets $\Gcal(c)$ for $L$ or for $\L$ coincide.

\lem{2.13} Let $\fun{\s}{\R}{\T}$ be $c$-minimal for $\L$; moreover, let $\s_t\in Mon$ for all $t\in\R$ and let
$$\liminf_{k\tends+\infty}\int_{-k}^k
[\L_c(t,\s_t,\dot\s_t)+\a(c)]\dt<+\infty  .  \eqno (2.18)$$
Let $\nu_0$ denote the Lebesgue measure on $I$. Then there is $q$, $c$-minimal for $L$, such that 
$(\pi\circ\s_t)_\sharp\nu_0=\d_{q(t)}$.

\proof We assert that it suffices to prove that
$$\int_\R\W(\s_t)\dt=0 . \eqno (2.19)$$
First of all, the integral above is well defined for any $\s_t$, though possibly equal to $-\infty$; this is because 
$\W\le 0$. Now, let $\s_t$ satisfy (2.19); since $\W\le 0$ and 
$\fun{}{t}{\W(\s_t)}$ is continuous (we recall that $\s_t$ is AC) we get that $\W(\s_t)=0$ for $t\in\R$. By our hypotheses on $W$, this implies that
$$\s_tx-\s_ty\in\Z\txt{for a. e.} (x,y)\in I\times I  .  $$
By Fubini and the fact that $\fun{}{t}{\s_tx}$ is AC for a. e. $x$, we can choose $x$ such that, setting, $q(t)=\s_tx$, we have

\noindent a) $\s_ty-q(t)\in\Z$ for a. e. $y\in I$.

\noindent b) $q$ is AC.

Now a) implies that $(\pi\circ\s_t)\nu_0=\d_{q(t)}$, which is part of the thesis; it remains to prove that $q$ is $c$-minimal.
To do this, we set $\tilde\s_t=q(t)$ and we see that, by point a) above, $\s_tx-\tilde\s_tx\in\Z$. We recall that, for a. e. $x\in I$, the function $\fun{}{t}{\tilde\s_tx-\s_tx}$ is $AC$; since it must take integer values, it is constant as a function of time. Thus, 
$\dot{\tilde\s}_tx=\dot\s_tx$ for a. e. 
$x\in I$; this implies the first equality below, while the second one comes from (2.19). 
$$\int_{-k}^k[\L_c(t,\s_t,\dot\s_t)+\a(c)]\dt=
\int_{-k}^k[\L_c(t,\tilde\s_t,\dot{\tilde\s}_t)+\a(c)]\dt=
\int_{-k}^k[L_c(t,q,\dot q)+\a(c)]\dt . $$
Now $\s$ is $c$-minimal by hypothesis; since by a)
$\tilde\s_{\pm k}-\s_{\pm k}\in L^2_\Z$, the first equality above implies that $\tilde\s$ is $c$-minimal too; this in turn implies, by the second equality above, that $q$ is $c$-minimal for the 
one-dimensional Lagrangian $L$.

We divide the proof of (2.19) into three steps. 

\noindent{\bf Step 1.} We begin with a much weaker fact than (2.19), i. e. that  
$$-\int_\R\W(\s_t)\dt<+\infty  .  \eqno (2.20)$$
Let $\s$ be as above; for starters, we prove that there is 
$C_1\ge 0$ such that
$$\int_{-k}^k[
\2 ||\dot\s_t||^2_{L^2(I)}-\inn{c}{\dot\s_t}_{L^2(I)}
-\V(\s_t)+\a(c)
]\dt\ge -C_1\quad\forall k\in\N   .   \eqno (2.21)$$
To prove (2.21), we recall a fact from [9]: there is $C_1\ge 0$ such that, for any absolutely continuous $q$ and $k\in\N$, 
$$\int_0^k[L_c(t,q,\dot q)+\a(c)]\dt\ge -C_1  .  \eqno (2.22)$$
To prove this, we recall that there is a $u\in C(S^1)$ which is a fixed point of $\hat\Lambda^k$ for all $k\in\N$; we defined the operator $\hat\Lambda$ before lemma 2.12. Let us suppose by contradiction that there is $q_k$ such that
$$\int_0^k[L_c(t,q_k,\dot q_k)+\a(c)]\dt\tends-\infty
\txt{as}k\tends+\infty . $$
Using $q_k$ as a test function in the definition of 
$\hat\Lambda^k(u)$, we see that $\hat\Lambda^k(u)$ is unbounded from below as $k\tends+\infty$; but this is impossible, since $\hat\Lambda^k(u)=u$ and $u$ is bounded.

Using Fubini, and the fact that the weak time derivative of $\s_t$ in 
$L^2$ coincides with $\frac{{\rm d}}{{\rm d}t}\s_tx$ for a. e. $x$, we see that
$$\int_{-k}^k[
\2 ||\dot\s_t||^2_{L^2(I)}-\inn{c}{\dot\s_t}_{L^2(I)}
-\V(\s_t)+\a(c)
]\dt=$$
$$\int_I\dx\int_{-k}^k[
\2 |\dot\s_tx|^2-c\cdot \dot\s_tx
-V(t,\s_tx)+\a(c)
]\dt\ge -C_1  $$
where the last inequality comes from (2.22). This proves (2.21).

Now, if (2.20) were false, (2.21) would imply 
$$\lim_{k\tends+\infty}\int_{-k}^k[
\L_c(t,\s_t,\dot{\s}_t)+\a(c)
]   \dt =  +\infty$$
contradicting (2.18).

\noindent{\bf Step 2.} We assert that (2.19) is true asymptotically, i. e. that
$$\lim_{|t|\tends+\infty}\W(\s_t)=0  .  \eqno (2.23)$$
To prove this, we recall that $\s_t$ is minimal; in particular, it satisfies $(ODE)_{Lag}$; thus,
$$||\ddot\s_t||_{L^2(I)}\le
||V^\prime(t,\s_t)||_{L^2(I)}+
||\int_I W^\prime(\s_t z-\s_t z^\prime){\rm d}z^\prime||_{L^2(I)}  .  
$$
Since $V^\prime$ and $W^\prime$ are bounded, we get that 
$$||\ddot\s_t||_{L^2(I)}\le C_2\qquad
\forall t\in\R  .  \eqno (2.24)$$
We assert that this implies that
$$||\dot\s_t||_{L^2(I)}\le C_3\qquad
\forall t\in\R    \eqno (2.25)$$
for a constant $C_3>0$.
Let us suppose by contradiction that 
$||\dot\s_{t_j}||_{L^2(I)}\tends+\infty$ for a sequence 
$t_j\tends+\infty$; we assume that $t_j\in[n_j,n_j+1]$. By (2.24), we have that 
$$\inf_{t\in[n_j,n_j+1]}||\dot\s_{t}||_{L^2(I)}\tends+\infty  
\txt{as} j\tends+\infty.  $$
Since $-W\ge 0$, we get the first inequality below; (2.21) implies the second one; the limit at the end follows by the formula above.
$$\int_{-k}^k
[\L_c(t,\s_t,\dot\s_t)+\a(c)]\dt\ge
\int_{-k}^{n_j}[\L_c(t,\s_t,\dot\s_t)+\a(c)]\dt+$$
$$\int_{n_j}^{n_j+1}\dt\int_I[L_c(t,\s_tz,\dot\s_tz)+\a(c)]\dr z+
\int^{k}_{n_j+1}[\L_c(t,\s_t,\dot\s_t)+\a(c)]\dt\ge$$
$$-2C_1+\int_{n_j}^{n_j+1}[
\2||\dot\s_t||_{L^2(I)}^2-\inn{c}{\dot\s_t}_{L^2(I)}-
\inn{1}{V(t,\s_t)}_{L^2(I)}+\a(c)
]\dt  \tends+\infty  .  $$
But this contradicts (2.18).

Since $W^\prime$ is bounded, we get that $\W^\prime$, the 
$L^2$-differential of $\W$, is bounded too; by (2.23) this implies that
$$\left\vert
\dif\W(\s_t)
\right\vert =
|\inn{\W^\prime(\s_t)}{\dot\s_t}_{L^2(I)}|
\le C_4\qquad\forall t\in\R  .  $$
Together with (2.20), this implies (2.23).

Now the idea is the following: if (2.19) did not hold, we could find a function $\tilde\s$ which coincides with $\s$ for $|t|\ge n$ and with the orbit of a single particle for $|t|\le n-1$; using step 2, we could prove that $\tilde\s$ has smaller action than $\s$, which is a contradiction since $\s$ is minimal. The next step tells us where to look for $\tilde\s$.

\noindent {\bf Step 3.} Let us suppose by contradiction that (2.19) does not hold; up to an integer translation in time, this means that 
$$-\int_0^1\W(\s_t)\dt\ge\e  \eqno (2.26)$$
for some $\e>0$. Once $\e$ is fixed in this way, we assert that there is $C>0$ for which the following holds: for each $\d>0$ we can find $k\in\N$, a set $A_k\subset I$, two functions 
$Z_{\pm k}\in L^2_\Z(I)$ and $z\in A_k$ such that
$$|(\s_{\pm k}x-Z_{\pm k}x)-(\s_{\pm k}z-Z_{\pm k}z)|\le 2C\sqrt\d \txt{if $k$ is large and $x\in A_k$.}   \leqno a)$$
$$|(\s_{\pm k}x-Z_{\pm k}x)-(\s_{\pm k}z-Z_{\pm k}z)|\le 4 \txt{if} 
x\not\in A_k  \leqno b)$$
and 
$$\int_{-k}^{k}L_c(t,\s_tz,\dot \s_tz)\dt\le
\int_{-k}^{k}[
\2||\dot\s_t||^2_{L^2(I)}-
\inn{c}{\dot\s_t}_{L^2(I)}-\V(t,\s_t)
]\dt+\frac{\e}{4}  .  \eqno (2.27)$$

We prove this fact. 
Chebishev's inequality and (2.23) imply this: let $\d,\chi>0$; then,  if $|t|$ is large enough, we have $-W(\s_tx-\s_tx^\prime)\le\d$ for 
$(x,x^\prime)\in (I\times I)\setminus B_{\d,\chi,|t|}$, with 
$\nu_0(B_{\d,\chi,|t|})<\chi$. By our hypotheses on $W$, this implies that 
$$|\s_tx-\s_tx^\prime|_{S^1}\le C\sqrt{\d}
\txt{for} 
(x,x^\prime)\in (I\times I)\setminus B_{\d,\chi,|t|}  .  \eqno (2.28)$$ 
We assert that this implies the following. Let 
$\d,\chi>0$; then for all $k\in\N$ large enough, we can find 
$q(\pm k)\in[\s_{\pm k}0,\s_{\pm k}1]$ and $Z_{\pm k}\in L^2_\Z$ such that, setting
$$A_{k}=\{
x\in I\st |\s_{\pm k}x-(Z_{\pm k}x+q(\pm k))|<C\sqrt\d
\}    $$
we have
$$\left\{
\eqalign{
\s_{\pm k}x-(Z_{\pm k}x+q(\pm k))\in [-1,1] \quad\forall (t,x)&{}\cr
\nu_0(A_k)\ge 1-\chi
\txt{for}k\tends+\infty  .  &{}
}
\right. \eqno (2.29) $$
Indeed, by Fubini, we can choose $x^\prime$ such that
$$\nu_0(\{
x\st (x,x^\prime)\in B_{\d,\chi,k}
\})  <\chi$$
and set $q(\pm k)=\s_{\pm k}x^\prime$; we choose 
$Z_{\pm k} x\in\Z$ in such a way that
$$|
\s_{\pm k}x-(Z_{\pm k}x+q(\pm k))
|=
|
\s_{\pm k}x-q(\pm k)
|_{S^1}  .  $$
This implies by (2.28) that 
$\nu_0(A_k)\ge 1-\chi$; the first estimate of (2.29) follows from the fact that $\s_{\pm k}\in Mon$. Moreover, since 
$q(\pm k)=\s_{\pm k}x^\prime\in[\s_{\pm k}0,\s_{\pm k}1]$, we can take $Z_{\pm k}\in\{ -1,0,1  \}$. 

We let $C_1$ be as in (2.22); we take $\chi\le\frac{\e}{2}$ in the second one of (2.29); taking $|k|$ large enough, we get that
$$C_1\nu_0(I\setminus A_{k})<\frac{\e}{8} . $$
We shall feel free to reduce $\chi$ (and thus to increase $k$) in the course of the proof.
Fubini implies the equality below; (2.22) implies the first inequality; the second one follows by the formula above.
$$\int_{-k}^k[
\2\norm{\dot\s_t}^2_{L^2(I)}-\inn{c}{\dot\s_t}_{L^2(I)}-
\V(t,\s_t)
]\dt  =$$
$$\int_{A_{k}}\dx\int_{-k}^k[
\2\abs{\dot\s_tx}^2-{c}\cdot{\dot\s_tx}-V(t,\s_tx)
]\dt+
\int_{I\setminus A_{k}}\dx\int_{-k}^k[
\2\abs{\dot\s_tx}^2-{c}\cdot{\dot\s_tx}-V(t,\s_tx)
]\dt\ge$$
$$\int_{A_{k}}\dx\int_{-k}^k[
\2\abs{\dot\s_tx}^2-{c}\cdot{\dot\s_tx}-V(t,\s_tx)
]\dt
-C_1\nu_0(I\setminus A_{k})  \ge  $$
$$\int_{A_{k}}\dx\int_{-k}^k[
\2\abs{\dot\s_tx}^2-{c}\cdot{\dot\s_tx}-V(t,\s_tx)
]\dt
-\frac{\e}{8} . $$
Thus (possibly reducing $\chi$ and thus increasing $k$), there is at least one 
$z\in A_{k}$, which depends on $k$, such that (2.27) holds.

Now point $a)$ follows because $z\in A_k$, and thus
$|\s_{\pm k}z-(Z_{\pm k}z-q(\pm k)|<C\sqrt\d$; moreover, for all 
$x\in A_k$, $|\s_{\pm k}x-(Z_{\pm k}x-q(\pm k)|<C\sqrt\d$. Point $b)$ follows in the same way, using the first formula of (2.29).

\noindent {\bf End of the proof of (2.19).} Let $z\in A_k$ be fixed as above; we set
$$\tilde\s_tx=\left\{\matrix{
\s_tz-(t+k-1)(\s_{-k}x-Z_{-k}x+Z_{-k}z-\s_{-k}z)&-k\le t\le -k+1\cr
\s_tz &|t|\le k-1\cr
\s_tz+(t-k+1)(\s_kx-Z_{k}x+Z_kz-\s_kz) &k-1\le t\le k  . 
}    \right.   $$
In other words, on $[-k,k]$ we squeeze together all the particles on 
$\s_tz$; this sends to zero the $\W$ term in the Lagrangian. We want to prove that, if (2.26) holds, then the action of 
$\tilde\s_t$ is smaller that the action of $\s_t$, contradicting the minimality of the latter.

We see that $\tilde\s_{\pm k}-\s_{\pm k}\in L^2_\Z$. 

Now points $a)$, $b)$ and the definition of $\tilde\s$ imply that
$$\left\{
\eqalign{
{}&
|\tilde\s_tx-\s_tz|\le 2C\sqrt\d 
\txt{if}
t\in[-k,-k+1], \quad 
x\in A_k \txt{and} k \txt{is large enough,}  \cr
{}&
|\tilde\s_tx-\s_tz|\le 2C\sqrt\d  
\txt{if} t\in[k-1,k],\quad 
x\in A_k \txt{and} k \txt{is large enough  .}   
}   
\right.    \leqno a^\prime)$$
and
$$\left\{
\eqalign{
{}&|\tilde\s_tx-\s_tz|\le 4
\txt{
if $t\in[-k,-k+1]$ and $k$ is large enough
}    \cr
{}&|\tilde\s_tx-\s_tz|\le 4   
\txt{
if $t\in[k-1,k]$ and $k$ is large enough.
}
}
\right.    \leqno b^\prime)   $$

By $a^\prime)$ and $b^\prime)$ we easily deduce that, possibly reducing again $\d$ and $\chi$, and then choosing a $k$ for which (2.29) holds, we have
$$\int_{[-k,-k+1]\cup[k-1,k]}|
\W(\tilde\s_t)-\W(\s_tz)
| \dt\le\frac{\e}{8},
\qquad
\int_{[-k,-k+1]\cup[k-1,k]}|
\V(t,\tilde\s_t)-\V(t,\s_tz)
| \dt\le\frac{\e}{8}    $$
and
$$\int_{[-k,-k+1]\cup[k-1,k]}[
\2\norm{\dot{\tilde\s}_t}^2_{L^2(I)}-\inn{c}{\dot{\tilde\s}_t}_{L^2(I)}-
(
\2\norm{\dot{\s}_tz}^2_{L^2(I)}-\inn{c}{\dot\s_tz}_{L^2(I)}
)
]\dt    \le\frac{\e}{8}  .  $$
Note that, as usual, we have identified $\s_tz\in\R$ with the function in $L^2(I)$ constantly equal to $\s_tz$; this, and the fact that $\tilde\s_tx\equiv\s_tz$ for $t\in[-k+1,k-1]$, imply the equality below. The first inequality below follows from the formula above and from the fact that
$$\int_{-k}^k\W(\s_tz)\dt=0$$
since $\fun{}{t}{\s_tz}$ is the motion of a single particle. The second inequality follows from (2.27); the last one follows from (2.26).
$$\int_{-k}^{k}[
\2||\dot{\tilde\s}_t||_{L^2(I)}^2-\inn{c}{\dot{\tilde\s}_t}_{L^2(I)}-\V(\tilde{\s}_t)-\W(\tilde{\s}_t)
]\dt=$$
$$\int_{-k}^k\L_c(t,\s_tz,\dot\s_tz)\dt+
\int_{[-k,-k+1]\cup[k-1,k]}\L_c(t,\tilde\s_t,\dot{\tilde\s}_t)\dt-
\int_{[-k,-k+1]\cup[k-1,k]}\L_c(t,\s_tz,\dot{\s}_tz)\dt\le$$
$$\int_{-k}^k L_c(t,\s_tz,\dot\s_tz)\dt+\frac{3}{8}\e\le
\int_{-k}^{k}[
\2||\dot\s_t||^2_{L^2(I)}-
\inn{c}{\dot\s_t}_{L^2(I)}-\V(t,\s_t)
]\dt+\frac{5}{8}\e\le
$$
$$\int_{-k}^k\L_c(t,\s_t,\dot\s_t)\dt-\frac{3}{8}\e  .  $$
Since $\tilde\s_{\pm k}-\s_{\pm k}\in L^2_\Z$, we have contradicted the minimality of $\s$.

\fin

We now give the proof of theorem 1, together with a more precise statement.

\thm{1} For any initial condition $M\in Mon$ and $c\in A$, there is $\s_t^c$, $c$-minimal for $\L$, such that $\s_0^c=M$, 
$\s_t^c\in Mon$ for $t\in[0,+\infty)$ and 
$$\lim_{t\tends+\infty}
\inf \{
dist_{weak}(\s_t^c,q(t))+||\dot\s_t^c-\dot q(t)||_{L^2(I)}\st
q\in\Gcal(c)  
\}    = 0  .  \eqno (2.30) $$
The limit above is uniform in the following sense: let $A\subset\R$ be compact and let $\e>0$. Then there is $T>0$ such that, for all $c\in A$ and all initial conditions $M\in Mon$,
$$\inf \{
dist_{weak}(\s_t^c,q(t))+||\dot\s_t^c-\dot q(t)||_{L^2(I)}\st
q\in\Gcal(c)  
\}<\e
\txt{if}t\ge T  .  \eqno (2.31)$$

\proof Let $U$ be a fixed point of $\Lambda_{c,\a(c)}$; we know that the $\inf$ in the definition of $\Lambda_{c,\a(c)}(U)$ is attained on some $\s^c_t$; we forgo the proof (which is identical to [4]) that $\s^c_t$ is defined for $t\in[0,+\infty)$ and that, for all 
$l\in\N$,
$$\Lambda^l_{c,\a(c)}(U)(M)=
\int_0^l[\L_c(t,\s^c_t,\dot\s^c_t)+\a(c)]\dt+U(\s^c_l) . 
\eqno (2.32)$$
The last formula easily implies that $\s_t^c$ is $c$-minimal on 
$[0,+\infty)$. By proposition 2.2, we can suppose that 
$\s_t^c\in Mon$ for $t\in[0,+\infty)$. It suffices to prove that 
$\s_t^c$ satisfies (2.30), or the stronger (2.31). 

Let us suppose by contradiction that (2.31) is false; then there is 
$\e>0$ and three sequences $n_k\nearrow+\infty$, 
$c_{k}\in A$ and $M_k\in Mon$ such that
$$\inf_{q\in\Gcal(c_{k})}\{
dist_{weak}(\s_{n_k}^{c_{k}}-q(n_k))+
||\dot\s_{n_k}^{c_{k}}-\dot q(n_k)||_{L^2(I)}
\}   \ge\e  .  \eqno (2.33)$$
Up to subsequences, we can suppose that $c_k\tends c\in A$. 
Since $\s^{c_k}_{n_k}\in Mon$, 
after an integer translation we can suppose that 
$\s^{c_k}_{n_k}$ is bounded. Moreover, $\s^{c_k}_{t+n_k}$ is uniformly H\"older from 
$[-n_k,+\infty)$ to $L^2(I)$; we saw this while proving proposition 2.2. Since $\s^{c_k}_{t+n_k}$ belongs to the locally compact space $Mon$, we can apply Ascoli-Arzel\`a\ as in the proof of proposition 2.2; we get that, up to subsequences, there is 
$\fun{\s}{\R}{L^2(I)}$ such that 
$\s^{c_k}_{t+n_k}\tends\s_t$ uniformly on the compact sets of 
$\R$. 

The orbit $\s$ is defined on $\R$ because each 
$\s^{c_k}_{t+n_k}$ is defined on $[-n_k,+\infty)$ and we are supposing that 
$n_k\nearrow+\infty$. Since $\s$ is the limit of a sequence of 
$c_{k}$-minimal orbits, it is easy to see that it is $c$-minimal. 
In particular, $\s$ satisfies $(ODE)_{Lag}$; thus, from the locally uniform convergence of $\s^{c_k}_{t+n_k}$ we deduce the locally uniform convergence of $\ddot\s^{c_k}_{t+n_k}$; from this, it follows through integration that $\dot\s^{c_k}_{t+n_k}$ converges, uniformly on compact sets, to $\dot\s_t$.

In lemma 2.8, it is easy to see that the Lipschitz constant $L(c)$ of 
$\Lambda_{c,\a(c)}U$ depends only on $c$ and is bounded on bounded sets. Since $\{ c_k \}$ converges, we can suppose that the fixed points of $\Lambda_{c_k,\a(c_k)}$ are $L$-Lipschitz for the same $L$; by (2.32), this implies that
$$\int_0^l[\L_{c_k}(t,\s^{c_k}_t,\dot\s^{c_k}_t)+\a(c)]\dt$$
is bounded independently on $l$ and $k$. Translating in time, it means that for any $a<b\in\Z$, for $k$ large we have that
$$\int_a^b[\L_{c_k}(t,\s^{c_k}_{t+n_k},\dot\s^{c_k}_{t+n_k})+\a(c)]\dt$$
is defined and bounded independently on $k$ and $b-a$. Taking limits under the integral sign, this implies that
$$\int_a^b[\L_c(t,\s_t,\dot\s_t)+\a(c)]\dt$$
is bounded too, i. e. that $\s$ satisfies (2.18). As a consequence, lemma 2.13 holds and 
$\s_t=q(t)$ with $q$ $c$-minimal. However, by (2.33) and uniform convergence of $(\s^{c_k}_{t+n_k},\dot\s^{c_k}_{t+n_k})$, we have that
$$\inf_{q\in\Gcal(c)}[
dist_{weak}(\s_0,q(0))+||\dot\s_0-\dot q(0)||_{L^2(I)}
]\ge\e$$
a contradiction.

\fin

\rem{2.15} We note a last connection with Aubry-Mather theory: if $W=0$ (i. e. if we are in the situation of [8]), and if $\s_t$ minimizes $\int_0^1\L_c$ among all curves such that 
$dist_{weak}(\s_1,\s_0)=1$, then 
$\dt\otimes(\s_t)_\sharp(\nu_0)$, a measure on 
$S^1\times S^1$, is a Mather $c$-minimal measure as defined in section 1. We leave the easy proof to the reader; other results along these lines are in [3].

\vskip 2pc
\centerline{\bf \S 3}
\centerline{\bf Birkhoff regions of instability}
\vskip 1pc

We prove theorem 2; we shall adapt the method of [2].

\vskip 1pc

\noindent{\bf Definitions.} Let $U$ be an open set of 
$S^1\times S^1$; we shall denote by $U$ also the open set of the points $(t,x)\in\R\times S^1$ such that $(t\mod 1,x)\in U$. Let 
$k\in\Z$. A one-form $\o$ on $\R\times S^1$ is called a 
$(U,k)$-step form if there is a closed form $\bar\o$ on 
$S^1\times S^1$ such that the restriction of $\o$ to $\{ t\le k-1 \}$ is zero, the restriction of $\o$ to $\{ t\ge k \}$ is $\bar\o$ and the restriction of $\o$ to the set 
$U\cup\{ t\le k-1 \}\cup\{ t\ge k \}$ is closed.

We shall also set
$$\tilde\Gcal(c)=\{
(t,q(t))\st q\in\Gcal(c),\quad t\in\R
\},\qquad
\tilde{\rm Lim}(c)=\{
(t,q(t))\st q\in{\rm  Lim}(c),\quad t\in\R
\}$$
where the sets $\Gcal(c)$ and ${\rm Lim}(c)$ have been defined in the introduction. 

\vskip 1pc

We omit the proof of the next lemma, which is a merger of point 2.6 and of theorem 6.3 of [2]; as usual, we denote by $[\o]$ the cohomology class of a closed form on $S^1\times S^1$.

\lem{3.1} Let us suppose that the interval $J\subset\R$ is a Birkhoff region of instability, and let $c\in J$. Then there is a neighbourhood $U$ of $\tilde\Gcal(c)$ such that, for any 
$k\in\Z$, there is a $(U,k)$-step form $\o$ with $[\bar\o]=[\dx]$.

\rm
\vskip 1pc

\noindent{\bf Definition.} We shall call adapted a neighbourhood $U$ of $\tilde\Gcal(c)$ as that of lemma 3.1.

\vskip 1pc

\lem{3.2} Let $\e\in(0,\frac{1}{8}]$, let the map $\fun{\s}{[a,b]}{\T}$ be continuous and let 
$\mu_t=(\pi\circ\s_t)_\sharp\nu_0$. Let us suppose that, for each 
$t\in[a,b]$, there is 
$x_t\in S^1$ such that 
$\mu_t(x_t-\e,x_t+\e)\ge 1-\e$. Then there is a function 
$\fun{q}{[a,b]}{\R}$ of class $C^1$ such that 

\noindent 1) $dist_{weak}(q(t),\s_t)<4\e$ and

\noindent 2) $|\pi q(t)-x_t|_{S^1}<4\e$.

Moreover, if $\s_b=a$ on $\T$, then we can choose $q$ in such a way that $q(b)=a$ on $S^1$.

\proof Since $\fun{\s}{[a,b]}{\T}$ is continuous, we get, as an easy consequence of dominated convergence, that the map 
$\fun{}{t}{\mu_t}$ is continuous for the weak-$\ast$ topology of measures. Let now $I\subset S^1$ be a fixed open arc; we recall that 
$$\mu_t(I)=\sup\int_{S^1}\phi\dr\mu_t$$
where the $\sup$ is taken among the continuous functions 
$\fun{\phi}{S^1}{[0,1]}$ which are zero outside $I$; this means that the map $\fun{}{t}{\mu_t(I)}$ is l. s. c., being the $\sup$ of a family of continuous functions. Thus, if 
$\mu_{t_0}(x_{t_0}-\e,x_{t_0}+\e)\ge 1-\e$, for $t_1$ close to 
$t_0$ we have that $\mu_{t_1}(x_{t_0}-\e,x_{t_0}+\e)\ge 1-2\e$.

By compactness, we can cover $[a,b]$ with a finite number of open intervals 
$I_i=(t_i-a_i,t_i+a_i)$, and we can find points $x_i\in S^1$ such that, if $t\in I_i$, then 
$\mu_{t}(x_{i}-\e,x_{i}+\e)\ge 1-2\e$. By the usual lemmas on coverings, we can suppose that

\noindent $a$) $t_{i}<t_{i+1}$.

\noindent $b$) $I_i\cap I_{i+1}\not =\emptyset$ and 
$I_i^c\setminus I_{i+1}\not =\emptyset\not=I_{i+1}\setminus I_i^c$.

\noindent $c$) $I_i\cap I_{i+2}=\emptyset$.

Let $t\in I_i\cap I_{i+1}$; such a point exists by $b$). We have that 
$\mu_t(x_i-\e,x_i+\e)\ge 1-2\e$ and 
$\mu_t(x_{i+1}-\e,x_{i+1}+\e)\ge 1-2\e$; since $\e\in(0,\frac{1}{8})$ and the total measure is $1$, the two intervals must intersect, and we get $|x_{i+1}-x_i|_{S^1}<2\e$. Thus, if 
$x\in[x_i,x_{i+1}]$, we have that
$$(x-3\e,x+3\e)\supset(x_{i}-\e,x_{i}+\e)
\txt{and}
(x-3\e,x+3\e)\supset(x_{i+1}-\e,x_{i+1}+\e)  .  $$
If $t\in I_i$ (or if $t\in I_{i+1}$), we already have an estimate on the measure of the two intervals on the right; thus,
$\mu_t(x-3\e,x+3\e)\ge 1-2\e$  if $t\in I_i\cup I_{i+1}$ and 
$x\in[x_i,x_{i+1}]$.

Now we take $\tilde q(t)\in C^1([a,b],S^1)$ such that 
$\tilde q(t_i)=x_i$ for all $i$ and 
$\tilde q(t)\in(x_i,x_{i+1})$ for $t\in(t_i,t_{i+1})$. By the last paragraph, for 
$t\in[a,b]$, $\mu_{t}(\tilde q(t)-3\e,\tilde q(t)+3\e)\ge 1-2\e$. We take a lift 
$q(t)$ of $\tilde q(t)$ to $\R$: we want to show that this function satisfies the thesis.

The equality below is proven in lemma 2.14 of [5]:
$$dist_{weak}(\s_t,q(t))=
\min\left\{
\int_{S^1\times S^1}|x-y|_{S^1}^2\dr\l(x,y)
\right\}  $$
where the minimum is taken on all the measures $\l$ on 
$S^1\times S^1$ which have $\mu_t$ as the first marginal and 
$\d_{\tilde q(t)}$ as the second one. If we take 
$\l=\mu_t\otimes\d_{\tilde q(t)}$, we get that
$$dist_{weak}(\s_t,q(t))\le
\int_{ ( \tilde q(t)-3\e,\tilde q(t)+3\e )\times S^1 }
|x-y|_{S^1}^2\dr\l(x,y)+
\int_{ ( \tilde q(t)-3\e,\tilde q(t)+3\e )^c\times S^1 }
|x-y|_{S^1}^2\dr\l(x,y)\le$$
$$9\e^2\mu_t(\tilde q(t)-3\e,\tilde q(t)+3\e)+
\mu_t((\tilde q(t)-3\e,\tilde q(t)+3\e)^c)\le
9\e^2+2\e . $$
Since $\e\le\frac{1}{8}$, we have that $9\e^2+2\e<4\e$ and point 1) follows.

We now recall that $\mu_t(\tilde q(t)-3\e,\tilde q(t)+3\e)\ge 1-2\e$ and 
$\mu_t(x_t-\e,x_t+\e)\ge 1-\e$; this implies as above that the two intervals must intersect, and thus point 2) follows.

As for the last assertion, we note that, if $\s_b=a$ on $\T$, then 
$\mu_b=\d_a$; we add to our covering of $[a,b]$ the interval 
$[b-\g,b+\g]$ on which $\mu_t(a-\e,a+\e)\ge 1-2\e$  and take the function $\tilde q$ in such a way that $\tilde q(b)=a$; then we are done.

\fin

\vskip 1pc
\noindent{\bf Definitions.} 

\noindent $\bullet$ {\bf The extension.} We shall say that 
$\fun{\s}{[a,b]}{\T}$ is $\e$-concentrated if it satisfies the hypotheses of lemma 3.2. If $\s_t$ is $\e$-concentrated and $q(t)$ is as in lemma 3.2, we shall say that the couple $(\s_t,q(t))$ is an extension of $\s_t$.

\noindent $\bullet$ {\bf The gap-filler.} Let $\{ c_i \}_{i\ge 0}$ be as in theorem 2. We fix $k_i\ge 0$ and we find $\{ d^i_j \}_{j=0}^{k_i-1}\subset[c_i,c_{i+1}]$ (or 
$\{ d^i_j \}_{j=0}^{k_i-1}\subset[c_{i+1},c_i]$ if $c_{i+1}\le c_i$) in such a way that 

\noindent $1)$ $d^i_{0}=c_{i}$; moreover, $\{ d^i_j \}_j$ is increasing if $c_i<c_{i+1}$, and decreasing if 
$c_i>c_{i+1}$.

$$|d^i_{j+1}-d^i_j|\le\frac{|c_{i+1}-c_i|}{k_i}
\txt{and} 
|d^i_{k_i-1}-c_{i+1}|\le\frac{|c_{i+1}-c_i|}{k_i}  .  \leqno 2)$$
We group the $\{ d^i_j \}_{i,j}$ into a unique sequence 
$\{ f_s \}_{s=0}^{+\infty}$ with the natural order; in other words,
$$f_0=c_0=d^0_0,\quad f_1=d^0_1,\quad f_2=d^0_2,\dots,\quad
f_{k_0-1}=d^0_{k_0-1},\quad f_{k_0}=c_1=d^1_0,\quad
f_{k_0+1}=d^1_1,\dots $$
We call such a sequence $\{ f_s \}_{s=0}^{+\infty}$ a 
$\{ k_i \}$-gap-filler of $\{ c_i \}$.

\vskip 1pc

\noindent $\bullet$ {\bf The form.} Let now  $\{ f_s \}_{s\ge 0}\subset J$ be a $\{ k_i \}$-gap-filler of 
$\{ c_i \}$ and let $\{ T_s \}_{s\ge 0}$ be a sequence of positive integers; we set $T_{-1}=0$ and, for $s\ge -1$,
$$\tilde T_s=\sum_{l=-1}^s T_l  .  $$
Using lemma 3.1, for $s\ge 0$ we can find adapted neighbourhoods $U_s$ of $\tilde\G(c_s)$ and 
$(U_s,\tilde T_s)$-step forms $\o_s$ such that
$$\o_s=\left\{\eqalign{
0&\txt{if}t\le\tilde T_s-1\cr
\bar\o_s&\txt{if}t\ge\tilde T_s
\txt{with}[\bar\o_s]=(f_{s+1}-f_s)[\dx]  .
}  \right.   \eqno (3.1)$$
We set $\o_{-1}=0$.

We define
$$\o=\sum_{l=0}^k\o_l$$
and we note that the form $c_0+\o$ 
satisfies $[c_0+\o]=f_s[\dx]$ for $s\ge 0$ and
$$t\in\left(
\tilde T_{s-1},\tilde T_s-1
\right)  .  $$
We refer the reader to the diagram below: on the middle line there are the times, on the upper ones the values of forms $\o_s$, and on the lower one the values of the form $c_0+\o$.
$$\matrix{
&{}&\omega_{-1}=0&{}&\omega_0=f_1-f_0&{}&\omega_1=f_2-f_1&{}&\omega_2=f_3-f_2&{}\cr
&\tilde T_{-1}=0&{}&\tilde T_0&{}&\tilde T_1&{}&\tilde T_2&{}&\tilde T_3\cr
&{} &c_0+\omega=f_0=c_0&{}&c_0+\omega=f_1&{}&c_0+\omega=f_2&{}&c_0+\omega=f_3&{}
}$$

\noindent $\bullet$ {\bf The set ${\cal D}$ of extensions.} Let 
$V_s$ be an open set such that 
${\cal G}(c_s)\subset V_s\subset\subset U_s$. Let
$$\d_s=\frac{1}{64}
\min\{
|x-y|_{S^1}\st (t,x)\in\bar V_s,\quad(t,y)\not\in U_s
\}      .  $$
We note that $\d_s\le\frac{1}{64}$, since the diameter of $S^1$ is smaller than $1$.
We shall consider the paths $\fun{\s}{[0,+\infty)}{\T}$ such that, setting as before 
$\mu_t=(\pi\circ\s_t)_\sharp\nu_0$,

\vskip 1pc

\noindent $A$) For all $s\ge 0$ and all 
$t\in[\tilde T_s-1,\tilde T_{s+1}-1]$, there is $x_t\in S^1$ such that 
$\mu_t([x_t-\d_s,x_t+\d_s])\ge 1-\d_s$. 

\vskip 1pc

\noindent $B$) $\m_t(\{ x\st(t,x)\in\bar V_s \})\ge 1-\d_s$ if $t\in[\tilde T_s-1,\tilde T_{s}]$ and $s\ge 0$. 

\vskip 1pc

Before ending the definition, we comment on conditions $A$) and $B$).

First of all, by point $A$), for $t\in[\tilde T_s-1,\tilde T_{s}]$ there is 
$x_t\in S^1$ such that $\mu_t([x_t-\d_s,x_t+\d_s])\ge 1-\d_s$; this and point $B$) imply that $[x_t-\d_s,x_t+\d_s]$ and 
$\{ x\st(t,x)\in\bar V_s \}$ must intersect; in other words, 
$x_t$ is in a $\d_s$-neighbourhood of $\{ x\st(t,x)\in\bar V_s \}$. 

Since $[x_t-\d_s,x_t+\d_s]\subset(x_t-2\d_s,x_t+2\d_s)$, we have by $A$) that 
$\s$ is $2\d_s$ concentrated; since $2\d_s\le\frac{1}{8}$, we can apply point 2) of lemma 3.2 with $\e=2\d_s$. We get that, for 
$t\ge\tilde T_0-1$, there is an extension 
$(\s,q)$ of $\s$ with $|q(t)-x_t|\le 8\d_s$; by the last paragraph, we get that $q(t)$ is in a $9\d_s$-neighbourhood of 
$\{ x\st(t,x)\in\bar V_s \}$; by the definition of $\d_s$, 
$(t,q(t))\in U_s$  for $t\in[\tilde T_s-1,\tilde T_s]$. 

We say that the couple $(\s,q)$, with $\s$ defined for $t\ge 0$ and $q$ defined for $t\ge\tilde T_0-1$, belongs to the set ${\cal D}$ if 
$\s$ satisfies $A$) and $B$), and if $q$ is an extension of $\s$ satisfying $|q(t)-x_t|_{S^1}\le 9\d_s$ for 
$t\in[\tilde T_s-1,\tilde T_{s+1}-1]$.

\lem{3.3} Conditions $A$) and $B$) are closed for the pointwise convergence of $\s_t$. In other words, if $\s^n_t$ satisfies $A$) and $B$) for $t\in[\tilde T_{s-1}-1,\tilde T_s-1]$, and if 
$\s^n_t\tends\s_t$ in $L^2$ for $t\in[\tilde T_{s-1}-1,\tilde T_s-1]$, then $\s_t$ satisfies $A$) and $B$) for $t\in[\tilde T_{s-1}-1,\tilde T_s-1]$.

\proof We recall that the map $\fun{}{\mu}{\mu(a,b)}$ is lower 
semi continuous for the weak$\ast$ topology of measures, and that the map $\fun{}{\s_t}{(\pi\circ\s_t)_\sharp\nu_0}$ is continuous from $L^2$ to the weak$\ast$ topology of measures. Let $\s^n_t$ and $\s_t$ be as in the statement of the lemma; let $\mu^n_t$ be the measure on 
$S^1$ induced by $\s^n_t$, and let $x^n_t\in S^1$ satisfy $A$). By what we just said, $\mu^n_t\tends\mu_t$ in the weak$\ast$ topology; up to taking a subsequence, we can suppose that 
$x^n_t$ tends to a point $x_t\in S^1$. We prove that $\mu_t$ and $x_t$ satisfy $A$). Let $\d^\prime>\d_s$; we have that
$$\mu_t([x_t-\d^\prime,x_t+\d^\prime]^c)\le
\liminf\mu^n_t([x_t-\d^\prime,x_t+\d^\prime]^c)\le
\liminf\mu_{t_n}([x_{t_n}-\d_s,x_{t_n}+\d_s]^c)\le\d_s  .  $$
The first inequality above comes from the fact that 
$[x_t-\d^\prime,x_t+\d^\prime]^c$ is an open arc of $S^1$ and the map 
$\fun{}{\mu}{\mu([x_t-\d^\prime,x_t+\d^\prime]^c)}$ is lower semicontinuous, as in lemma 3.2. The second inequality comes from the fact that 
$[x^n_{t_n}-\d_s,x^n_{t_n}+\d_s]\subset[x_{t}-\d^\prime,x_t+\d^\prime]$ for 
$n$ large enough, since $\d^\prime>\d_s$ and $x^n_t\tends x_t$; the third one, from the fact that $\s^n_t$ satisfies $A$) . Since 
$\d^\prime>\d_s$ is arbitrary, we get that $A$) holds for $\s_t$. 

The proof of $B$) is analogous.

\fin

Let now $\fun{\s}{[0,+\infty)}{\T}$ be an orbit which satisfies $B$); let it satisfy $A$) for $x_t\in S^1$, and also for $\tilde x_t\in S^1$. Let $(\s,q)$ be an extension with $|\pi\circ q(t)-x_t|_{S^1}<9\d_s$, and let $(\s,\tilde q)$ be another extension with 
$|\pi\circ\tilde q(t)-\tilde x_t|_{S^1}<9\d_s$. Since $A$) holds for $x_t$ and $\tilde x_t$, we get that $|x_t-\tilde x_t|_{S^1}\le 2\d_s$; in particular, $|\pi\circ\tilde q(t)-x_t|_{S^1}<11\d_s$. Since we also have that 
$|\pi\circ q(t)-x_t|_{S^1}<9\d_s$, we get that 
$|\pi\circ\tilde q(t)-\pi\circ q(t)|_{S^1}\le 20\d_s$ for 
$t\ge\tilde T_0-1$. Since $\d_s\le\frac{1}{64}$, we have that 
$20\d_s<\2$; i. e., if 
$\fun{\pi}{\R}{S^1}$, $q(t)$ and $\tilde q(t)$ are in a ball on which a determination of $\pi^{-1}$ is defined.

This has a standard consequence for the lifts to the universal cover: up to adding an integer to $\tilde q$, there is a unique lift of $x_t$ to $\R$ such that
$|q(t)-x_t|<11\d_s$ and 
$|\tilde q(t)-x_t|<11\d_s$ for all
$t\in[\tilde T_{s}-1,\tilde T_{s+1}-1]$ and all $s\in\N$.

\lem{3.4} Let $\s$ satisfy points $A$) and $B$) above, and let 
$(\s,q)$, $(\s,\tilde q)\in{\cal D}$ be two extensions of $\s$; up to adding an integer to $\tilde q$, we suppose that $q$ and 
$\tilde q$ satisfy the property above. Then, the curves $\fun{}{t}{(t,q(t))}$ and $\fun{}{t}{(t,\tilde q)}$ are homotopic in 
$${\cal U}\colon=
\{ \tilde T_{-1}\le t\le\tilde T_0-1 \}\cup U_0\cup
\{ \tilde T_{0}\le t\le\tilde T_1-1 \}\cup U_1\cup\dots\subset
\R\times\R^n.  $$

\proof Before stating the lemma, we saw that $x_t$ can be lifted in a unique way to $\R$ in such a way that, for 
$t\in[\tilde T_s-1,\tilde T_s]$, we have  $|q(t)-x_t|\le 11\d_s$ and 
$|\tilde q(t)-x_t|\le 11\d_s$. Since $x_t$ is in a 
$\d_s$-neighbourhood of $\{ x\st(t,x)\in\bar V_s \}$ (we noticed this before while defining {\cal D}), we get that, for $\l\in[0,1]$, 
$(1-\l)q(t)+\l\tilde q(t)$ is in a $12\d_s$-neighbourhood of 
$\{ x\st (t,x)\in\bar V_s \}$; by the definition of $\d_s$, this implies that $(t,(1-\l)q(t)+\l\tilde q(t))\in U_s$ for $\l\in[0,1]$, and this ends the proof.

\fin

\noindent{\bf Definition.} We shall say that $(\s,q)\in{\cal D}$ is a minimizer of $\L_{c_0}-\o$ if, for all $t_1>T_0$, the following happens. Let 
$(\tilde\s,\tilde q)\in{\cal D}$ satisfy  
$\tilde\s_{0}-\s_{0}\in L^2_\Z$, 
$\tilde\s_{t_1}-\s_{t_1}\in L^2_\Z$, and
$\tilde q(t_1)=q(t_1)$. Then,
$$\int_{0}^{t_1}
\L_{c_0}(t,\s_t,\dot\s_t)\dt  -
\int_{\tilde T_0-1}^{t_1}
\o(t,q(t))\cdot(1,\dot q(t))
\dt \le
\int_{0}^{t_1}
\L_{c_0}(t,\tilde\s_t,\dot{\tilde\s}_t)\dt-
\int_{\tilde T_0-1}^{t_1}
\o(t,\tilde q(t))\cdot
(1,\dot{\tilde q}(t))
\dt   .   $$

\vskip 1pc

\lem{3.5} Let $M\in Mon$. Then there is 
$(\s,q)\in{\cal D}$ such that

\noindent 1) $\s_0=M$.

\noindent 2) $\s_t\in Mon$ for all $t\in[0,+\infty)$.

\noindent 3) $(\s,q)$ minimizes $\L_{c_0}-\o$ in ${\cal D}$.

\proof We begin to prove this: for $\tilde T_s$ defined as above, there is $(\s^s,q^s)\in{\cal D}$ which minimizes 
$\L_c-\o$ on $[0,\tilde T_s]$ with boundary conditions 
$\s^s_0=M$, $\s^s_{\tilde T_s}=0$ and $q^s(\tilde T_s)=0$. Naturally, when we say that $\s^s_0=M$ or $\s^s_{\tilde T_s}=0$, we mean equality on $\T$, i. e. up to adding an element of 
$L^2_\Z$; also $q^s(\tilde T_s)=0$ is an equality on $S^1$.

The proof of this is similar to the one of proposition 2.2. We begin to tackle the finite-dimensional problem.

Let $P_n$ be the projection of section 1; we want to connect $P_nM$ and $0$ with a minimal path $(\tilde\s,\tilde q)\in{\cal D}$, with ${\tilde\s_t}\in{\cal C}_n$. 
Let us set, as in section 1, $P_n\s_t=D_n(z_1(t),\dots,z_n(t))$; we are thus minimizing the functional
$$I(z)=\int_0^{\tilde T_s}[
\frac{1}{n}\sum_{i=1}^n |\dot z_i|^2-
\frac{c_0}{n}\sum_{i=1}^n\dot z_i-
\frac{1}{n}\sum_{i=1}^n V(t,z_i)-
\frac{1}{2n^2}\sum_{i,j=1}^nW(z_i-z_j)]\dt-
\int_{\tilde T_0-1}^{\tilde T_s}
\o(t,q(t))\cdot(1,\dot q(t))
\dt$$
over all couples $(z,q)$ belonging to ${\cal D}$. As for the boundary conditions, we ask that $z(0)=P_nM$ and 
$z(\tilde T_s)=0$ in $\T$; equivalently, we ask that, if 
$P_nM=(\bar z_1,\dots,\bar z_n)$, then
$\pi\circ z_i(0)=\pi(\bar z_i)$ and $\pi\circ z_i(\tilde T_s)=0$ for $i=1,\dots,n$. Moreover, we ask that $\pi\circ q(\tilde T_s)=0$. Such an extension exists because of the last assertion of lemma 3.2.

Note that we have called the functional $I(z)$ even though, nominally, it depends on $(z,q)$. However, by lemma 3.4, if 
$(z,q)$ and $(z,q_1)$ are two extensions of $z$, then $q$ and 
$q_1$ are homotopic in ${\cal U}$; since $\o(\tilde T_0-1,\cdot)=0$ and $q(\tilde T_s)=q_1(\tilde T_s)$, this implies that 
$$\int_{\tilde T_0-1}^{\tilde T_s}\o(t,q(t))\cdot(1,\dot q(t))\dt=
\int_{\tilde T_0-1}^{\tilde T_s}\o(t,q_1(t))\cdot(1,\dot q_1(t))\dt  .  $$
In other words, $I(z)$ does not depend on the particular extension $(z,q)$ of $z$ we choose: it is a function only of $z$.

By lemma 3.3, the fact that $z$ satisfies $A$) and $B$) is equivalent to the fact that $(t,z(t))\in K_1$, where $K_1$ is a closed set in $\R\times\R^n$. In other words, we are dealing with a "minimization with obstacle" problem: we are minimizing $I$ among all $z\in AC([0,\tilde T_s],\R^n)$ such that 
$(t,z(t))$ belongs to a closed set $K_1$. It is standard (see below for a proof) that such problems admit a minimum, provided they are coercive, and this is what we prove next.

By the definition of $\o$, we have that
$$
\int_{\tilde T_0-1}^{\tilde T_s}\o(t,q(t))\cdot(1,\dot q(t))\dt=
\int_{\tilde T_0-1}^{\tilde T_1-1}\o_0(t,q(t))\cdot(1,\dot q(t))\dt+
\int_{\tilde T_1-1}^{\tilde T_2-1}\o_1(t,q(t))\cdot(1,\dot q(t))\dt+\dots+
$$
$$
\int_{\tilde T_{s-1}-1}^{\tilde T_s-1}\o_{s-1}(t,q(t))\cdot(1,\dot q(t))\dt+
\int_{\tilde T_s-1}^{\tilde T_s}\o_s(t,q(t))\cdot(1,\dot q(t))\dt  .
$$
Since the forms $\o_1,\dots,\o_s$ are finite in number, and each of them is bounded, we get that  $I$ is coercive; actually, we get that, if $z^k$ is a minimizing sequence for $I$, then
$$\frac{1}{n}\int_0^{\tilde T_s}\sum_{i=1}^n|\dot z^k_i|^2\dt\le C_4
\eqno (3.2)$$
for some $C_4$ independent on $n$ and $k$. In particular, $z^k$ is uniformly H\"older; since $z^k(T_0)=P_nM$, we have that 
$z_k(0)$ is bounded and thus, that $z_k$ is bounded on 
$[0,\tilde T_s]$. By Ascoli-Arzel\`a, 
$z^k$ converges, up to subsequences, to a limit $\tilde z$; since $K_1$ is closed, we get that $(t,\tilde z(t))\in K_1$ for all 
$t\in[0,T]$. It is standard (see section 2 for the proof of a similar fact) that $I$ is lower semicontinuous under uniform convergence; by this and lemma 3.4, if we take any extension 
$(\tilde z,\tilde q)\in{\cal D}$ with 
$\pi\circ\tilde q(\tilde T_s)=0$, $(\tilde z,\tilde q)$ will minimize 
$I(z)$ among all couples $(z,q)\in{\cal D}$ such that 
$z\in AC([0,\tilde T_s],\R^n)$ and $(z,q)$ satisfies the boundary conditions above.

Now we want to prove the assertion we made at the beginning, i. e. the existence of a minimum when $\s$ takes values in $L^2(I)$, not in 
${\cal C}_n$. Let us call $(\tilde\s^n,\tilde q^n)$ the minimal couple we found above, and let us set 
$\tilde\s^n=D_n(\tilde z_1,\dots,\tilde z_n)$. We note that the argument of lemma 2.7 continues to hold: indeed, this argument consisted in rearranging the indices of $(\tilde z_1,\dots,\tilde z_n)$; but this has no effect on properties $A$) and $B$). Thus we can suppose that 
$\tilde\s^n_t\in Mon$ for 
$t\in[0,\tilde T_s]$; it follows by (3.2) that the $\2$-H\"older norm of $\tilde\s^n$ is bounded uniformly in $n$. This implies as in the proof of proposition 2.2 that, up to subsequences, 
$\tilde\s^n\tends\s^{\tilde T_s}$. By lemma 3.3, 
$A$) and $B$)  are closed conditions, and thus $\s^{\tilde T_s}$ continues to satisfy them. Again, we take any extension 
$(\s^{\tilde T_s},q^{\tilde T_s})\in{\cal D}$ with 
$\pi(q^{\tilde T_s}(\tilde T_s))=0$, and that will be minimal by the lower semicontinuity of $I$ and lemma 3.4.

We now note that, by an argument similar to the one above, 
$\s_t$ is bounded in the $C^{0,\2}_{loc}$ topology; in other words, on any fixed set $[0,T]$, the $\2$-H\"older norm of $\s^{T_s}$ is bounded in $s$; since 
$\s_t^{\tilde T_s}\in Mon$ for $t\in[0,\tilde T_s]$, we can use Ascoli-Arzel\`a\ as in the proof of proposition 2.2 and get that, up to subsequences, $\s^{\tilde T_s}\tends\s$ uniformly on compact sets. Now we see as before that $\s$ is minimal on $[0,+\infty)$, and we are done.

\fin

We omit the proof of the next lemma, since it is identical to that of theorem 1.

\lem{3.6} Let $\{ f_s \}_{s\ge 0}$ be as above and let 
$\{ \g_s \}_{s\ge 0}$ be a sequence in $(0,1)$. Then, if we choose the times $\{ T_s \}_{s\ge 0}$ large enough, the following happens. Let $\fun{\s}{[0,+\infty)}{\T}$ minimize $\L_{c_0}-\o$ in 
${\cal D}$; in view of lemma 3.5, we shall suppose that 
$\s_t\in Mon$ for $t\in[0,+\infty)$. Then, for any 
$t\in[\tilde T_{s-1},\tilde T_s]$, $s\ge 1$, we can find $\bar q$, minimal for $L_{c_0}-\o$, satisfying
$$dist_{weak}(\s_t, \bar q(t))+
\norm{\dot\s_t-\dot{\bar q}(t)}_{L^2(I)}
  <\frac{\g_s}{4}  .  $$

\rm

\vskip 1pc

\lem{3.7} There are $\Gamma_i>0$ such that the following holds. Let $\{ d^i_j \}$ be a $\{ k_i \}$-gap-filler of 
$\{ c_i \}$, let $k_i\ge\Gamma_i$ and let 
$T_0,\dots,T_{k_1}\ge\Gamma_1$, 
$T_{k_1+1},\dots,T_{k_2}\ge\Gamma_2$, etc...
Let $(\s,q)$ minimize 
$\L_c-\o$ in $\cal D$.
Then $\s_t$ is a solution of $(ODE)_{Lag}$.

\proof We know from lemma 3.6  that, if we choose $\Gamma_i$ large enough, for any $t\in[\tilde T_{s-1},\tilde T_s]$, $s\ge 1$ there is 
$\fun{\bar q}{[0,+\infty)}{S^1}$, minimal for $L_c-\o$, such that
$$dist_{weak}(\s_t,\bar q(t))<\frac{\g_s}{4} . $$
Now, $\bar q$ is a minimal orbit of the one-dimensional Lagrangian $L_c-\o$; for this Lagrangian it has been proven in [2] that, if $T_{s-1},T_s,T_{s+1}$ are large enough and 
$a\in [\tilde T_{s-1},\tilde T_{s}]$, $s\ge 1$,
there is $\tilde q\in{\cal G}(f_s)$ such that
$$\sup_{t\in[a,a+1]}|\bar q(t)-\tilde q(t)|_{S^1}<\frac{\g_s}{4}  .  $$
Together with the last formula, this implies that, if we choose 
$\Gamma_i$ large enough, we have that, for $s\ge 1$ and 
$a\in [\tilde T_{s-1},\tilde T_{s}]$, there is $\tilde q\in{\cal G}(f_s)$ such that 
$$\sup_{t\in[a,a+1]}dist_{weak}(\s_t,\tilde q(t))<\frac{\g_s}{2}  .  
\eqno (3.3)$$

Thus, if we choose $\g_s\le\d_s$, we have $A^\prime$) below; possibly reducing $\g_s$, we get $B^\prime$).

\vskip 1pc

\noindent $A^\prime$) $\s_t$ is $\frac{1}{2}\d_s$-concentrated for 
$t\in[\tilde T_s,\tilde T_{s+1}]$ and $s\ge 0$.

\vskip 1pc

\noindent $B^\prime$) 
$\mu_t(\bar V_s\cap(\{ t \}\times S^1))\ge 1-\2\d_s$ for
$t\in[\tilde T_s-1,\tilde T_s]$ and $s\ge 0$.

\vskip 1pc

In particular, $\s$ satisfies points $A$) and $B$) above; thus, we can find an extension $(\s,q)$ such that $(t,q(t))\in U_s$ for 
$t\in[\tilde T_s-1,\tilde T_s]$.

Now we can show that the Euler-Lagrange equation holds. Let 
$T>\tilde T_0$, let
$\psi\in C^1([0,T],L^2(I))$ and let us suppose that 
$\psi_{0}=\psi_{T}=0$. Let us set 
$\s^\l=\s+\l\psi$; the boundary conditions on $\psi$ imply that 
$$\s^\l_{0}=\s_{0} \txt{and}
\s^\l_{T}=\s_{T}  .  $$
Moreover, 
$$\sup_{t\in[0,T]}||
\s^\l_t-\s_t
||_{L^2(I)}     \tends 0   \txt{as}  \l\tends 0  .  $$
The formula above and $A^\prime$), $B^\prime$) imply that, if $\l$ is small enough, then $\s^\l$ satisfies points $A$) and $B$). We have seen that this implies that  $\s_t^\l$ has an extension 
$(\s_t^\l,q^\l(t))$ with 
$(t,q^\l(t))\in U_s$ for $t\in [\tilde T_s-1,\tilde T_s]$. Moreover, arguing as in lemma 3.2, we can require that $q^\l(T)=q(T)$. By lemma 3.4, $q$ and $q^\l$ are homotopic in ${\cal U}$; moreover, $\o(T_0-1,\cdot)=0$. Thus,
$$\int_{\tilde T_0-1}^{T}
\o(t,q(t))\cdot(1,\dot q(t))\dt=
\int_{\tilde T_0-1}^{T}\o(t,q^\l(t))\cdot(1,\dot q^\l(t))\dt  .  $$
This and the fact that $\s_t$ is minimal for $\L-\o$ imply that
$$\int_{0}^{T}
\L_{c_0}(t,\s_t^\l,\dot\s_t^\l) \dt-
\int_{0}^{T}\L_{c_0}(t,\s_t,\dot\s_t)\dt\ge 0 . $$
Differentiating in $\l$, we get that $\s$ solves $(ODE)_{Lag}$.

\fin

\lem{3.8} For any $c\in\R$ and any $\e>0$, there is 
$T_\e(c)\in\N$ with the following properties: if 
$\fun{\s}{[0,T_\e(c)]}{Mon}$ is $c$-minimal for $\L$, then there is a time $t\in[0,T_\e(c)]$ such that
$$\inf_{q\in{\cal M}(c)}[
dist_{weak}(\s_t,q(t))+
\norm{\dot\s_t-\dot q(t)}_{L^2(I)}
]\le\e  .  $$
Moreover, if $\fun{\s}{[0,+\infty]}{Mon}$ is $c$-minimal for 
$\L$, then 
$$\inf_{q\in{\rm Lim}(c)}[
dist_{weak}(\s_t,q(t))+
\norm{\dot\s_t-\dot q(t)}_{L^2(I)}
]\le\e  \txt{for}
t\ge \2 T_\e(c)   .  $$

\proof We begin with the first statement. By lemma 2.9 of [2], there is $\tilde T_\e(c)\in\N$ with the following property: if $T\in\N$, if 
$\fun{q}{[T,T+\tilde T_\e(c)]}{S^1}$ is $c$-minimal for $L$, then there is $t\in[T,T+\tilde T_\e(c)]$ and $\tilde q\in{\cal M}(c)$ such that
$$|q(t)-\tilde q(t)|+|\dot q(t)-\dot{\tilde q}(t)|<\frac{\e}{2}  .  $$
By theorem 1, for $\g>0$, there is $T\in\N$ and $q$ $c$-minimal such that
$$dist_{weak}(\s_T,q(T))+||
\dot\s_T-\dot q(T)
||_{L^2(I)}\le\g  .  $$
If we take $\g$ small enough and $T_\e(c)=T+\tilde T_\e(c)$, the thesis follows from the last two formulas and continuous dependence on the initial conditions.

The proof of the second statement is similar to the first one: $\s_t$ accumulates, for $t$ large, on a $c$-minimal orbit of $L$; this orbit accumulates on ${\rm Lim}(c)$ by definition of the latter. 

\fin

\vskip 1pc

\noindent{\bf Proof of theorem 2.} Let $\{ f_s \}$ be a 
$\{ k_i \}$-gap-filler of $\{ c_i \}$; let $\G_i$ satisfy the hypotheses of lemma 3.7, and let $k_i=\G_i$. We define the form $\o$ as above, choosing the times $T_s$ in the following way. When 
$f_s=d^i_j$ and $j\not=0$, we take $T_s=k_i$. When $f_s=d^i_0$, we take
$$T_s=t_i^\pprime-t_i^\prime-(k_i-1)k_i+
(t_{i+1}^\prime-t_i^\pprime) .  $$
Here, $t_i^\prime$ and $t_i^\pprime$ are as in the hypotheses of theorem 2; since we want that $T_s\ge k_i$, we ask that
$$t_{i+1}^\prime-t_i^\pprime\ge T(c_i,c_{i+1})
\colon= k_i^2  .  $$
Moreover, we ask that 
$$t_i^\pprime-t_i^\prime\ge T_{\e_i}(c_i)$$
where $T_{\e_i}(c_i)$ is defined as in lemma 3.8. 

Now lemma 3.5 holds for any choice of the $T_i$, and thus there is $(\s,q)$ minimal in ${\cal D}$ for $\L-\o$. We have chosen 
$k_i$ and $T_i$ in such a way that they satisfy the hypotheses of lemma 3.7. As a consequence, $\s$ solves 
$(ODE)_{Lag}$. Formula (1) of theorem 2 holds by our choice of 
$T_\e(c_i)$ and lemma 3.8. As for the last statement of theorem 2, it follows from the last statement of lemma 3.8.

\fin

\vskip 1pc
\noindent{\bf Acknowledgements.} We thank the referees for helpful comments and remarks.

\vskip 2pc
\centerline{\bf Bibliography}


\noindent [1] R. A. Adams, Sobolev spaces, Academic Press, 1975, New York.

\noindent [2] P. Bernard, Connecting orbits of time-dependent Lagrangian systems, Annales de l'Institut Fourier, Grenoble, 
{\bf 52}, 1533-1568, 2002.

\noindent [3] P. Bernard, B. Buffoni, Optimal mass transportation and Aubry-Mather theory, J. Eur. Math. Soc., {\bf 9}, 85-121, 2007.

\noindent [4] A. Fathi, Weak KAM theorem in Lagrangian dynamics, Fourth preliminary version, mimeographed notes, Lyon, 2003.

\noindent [5] W. Gangbo, A. Tudorascu, Lagrangian dynamics on an infinite-dimensional torus; a weak KAM theorem, Adv. Math., 
{\bf 224},  260-292, 2010.

\noindent [6] W. Gangbo, T. Nguyen, A. Tudorascu, Hamilton-Jacobi equations in the Wasserstein space, Methods Appl. Anal., 
{\bf 15}, 155-183, 2008.

\noindent [7] W. Gangbo, A. Tudorascu, A weak KAM theorem: from finite to infinite dimension, in Optimal transportation, Geometric and functional inequalities, L. Ambrosio editor, Pisa, 2010.


\noindent [8] J. N. Mather, Action minimizing invariant measures for 
positive-definite Lagrangian Systems, Math. Zeit., {\bf 207}, 
169-207, 1991.


\noindent [9] J. N. Mather, Variational construction of connecting 
orbits, Ann. Inst. Fourier, {\bf 43}, 1349-1386, 1993.

\end

\lem{3.6} Let $\d,\tilde T>0$ be as above and let $M\in Mon$. Then, the minimal orbit of lemma pippo is $\frac{1}{16}$-concentrated on $[T,+\infty)$.

\proof This follows from a proof.

\fin

\lem{3.7} Let $c_0,c_1\in J$. Then there exists an integer 
$T(c_0,c_1)$ and a one-form $\o$ on $\R\times S^1$ such that the two points below hold. 

\noindent 1) The restriction of $\o$ to $\{  t\le 0  \}$ is zero and the restriction of $\o$ to $\{ t\ge T(c_0,c_1) \}$ satisfies 
$[\bar\o]=(c_1-c_0)[\dx]$.

\noindent 2) If $\fun{\s}{[c,d]}{L^2(I)}$ is $\e_0$-concentrated and minimizes $\L-\o$, then $\s$ is an extremal of $\L$ if $[c,d]$ contains $[0,T(c_0,c_1)]$.

\proof Again, we only sketch the proof of this lemma, referring the reader to lemma 2.8 of [2] for details.

By compactness, it is possible to find a finite partition
$$c_0=\bar c_0<\bar c_1<\dots<\bar c_k=c_1$$
of $[c_0,c_1]$ in such a way that $|c_{i+1}-c_i|\le\g(c_i)$, where 
$\g(c_i)$ is as in lemma 3.3. By lemma 3.3, we can find three sequences, one of adapted neighbourhoods 
$U_i\supset\Gcal(\bar c_i)$, another one of $U_i$-step forms 
$\o_i$ and a third one of times $T_i$ such that
$$\o_i=\left\{\eqalign{
0&\txt{if}t\le\sum_{l=0}^iT_l\cr
\bar\o_i=&\txt{if}t\ge\sum_{l=0}^iT_l+1
\txt{with}[\bar\o_i]=(c_{i+1}-c_i)[\dx]  .
}  \right.   \eqno (3.1)$$
Moreover, we have that, if $\s$ satisfies the hypotheses of point 2), then $\s$ is an extremal of $\L$ on each $[T_i,T_{i+1}]$, and thus is an extremal of $\L$ on $[c,d]$.

We now set $T_{-1}=0$ and we note that the form
$$\o=c\dx+\sum_{l=1}^k\o_l$$
satisfies $[\o_i]=c_i[\dx]$ for
$$t\in\left(
\sum_{l=-1}^{i-1}T_l+1,\sum_{l=-1}^iT_l
\right)  .  $$
Assertion 1) follows by (3.1).

\fin

\lem{2.8} Let $M\in L^2(I)$. Let $q=(q_1,\dots,q_n)$ be minimal in the second formula of (2.6) (or, equivalently, let
$\g^n=D_nq$ be minimal in the first one); then there is $C>0$, independent on $n$, such that 
$$\int_0^1||\dot\g^n_t||^2_{L^2(I)}\dt=
\frac{1}{n}\int_0^1\sum_{i=1}^n|\dot q_i|^2\dt\le C  .  $$

\proof Let us consider a constant trajectory $\tilde q$:
$$D_n\tilde q(t)=P_nM\txt{for} t\in[0,1] . $$
Formula (2.6) implies the first inequality below.
$$\d_n\le
\int_0^1L_{n,c}(t,\tilde q,\dot{\tilde q})\dt+U(D_n\tilde q(1))=$$
$$-\frac{1}{n}\int_0^1\sum_{i=1}^nV(t,\tilde q_i)\dt-
\frac{1}{2n^2}\int_0^1\sum_{i,j=1}^nW(\tilde q_i-\tilde q_j)\dt+
U(D_n\tilde q(1))\le C_1  .  $$
The last inequality comes from the fact that $V$ and $W$ are bounded by hypothesis and $U$ is bounded by lemma 2.3. 

If $q$ is minimal in (2.6), the formula above implies
$$\frac{1}{n}\int_0^1\sum_{i=1}^n\left[
\2 |\dot q_i|^2-c\dot q_i
\right]  \dt-
\frac{1}{n}\int_0^1\sum_{i=1}^nV(t, q_i)\dt-
\frac{1}{2n^2}\int_0^1\sum_{i,j=1}^nW(q_i-q_j)\dt+
U(D_n\tilde q(1))\le C_1  .  $$
Using again the fact that $V$, $W$ and $U$ are bounded, we get that
$$\frac{1}{n}\int_0^1\sum_{i=1}^n\left[
\2 |\dot q_i|^2-c\dot q_i
\right]  \dt   \le   C_2    $$
or equivalently
$$\int_0^1[
\2 ||\dot\g^n_t||^2_{L^2(I)}-\inn{c}{\dot\g^n_t}_{L^2(I)}
]  \dt\le C_2  .  $$
Now the thesis follows by the H\"older inequality.

\fin

Since $z$ satisfies points $A$) and $B$) above, lemma 3.2 implies  that, on $S^1$, $(t,q(t))\in U_i$ if 
$t\in[\tilde T_{i}-1,\tilde T_i]$; since 
$\o_i$ closed on this set, the integral of 
$\o_i(t,q(t))\cdot(1,\dot q(t))$ on $[\tilde T_{i-1},\tilde T_i]$ depends only on the position of  $q(\tilde T_{i-1})$ and $q(\tilde T_i)$ on the universal cover $\R$; since the forms 
$\o_s$ are bounded, the last formula implies the bound
$$\left\vert
\int_{\tilde T_0-1}^{\tilde T_s}\o(t,q(t))\cdot(1,\dot q(t))\dt
\right\vert
\le C_1
(
|q(\tilde T_1)-q(\tilde T_0)|+|q(\tilde T_2)-q(\tilde T_1)|+\dots+
|q(\tilde T_s)-q(\tilde T_{s-1})|
)  .  $$
Since $(z,q)$ is an extension of $z$, the first inequality below follows easily; the second one is H\"older; we have denoted by 
$||\cdot||$ the Euclidean norm on $\R^n$.
$$|q(\tilde T_i)-q(\tilde T_{i-1})|\le
\frac{C_2}{n}\sum_{j=1}^n|z_j(\tilde T_i)-z_j(\tilde T_{i-1})|\le
C_2[\tilde T_i-\tilde T_{i+1}]^\2\cdot
\left[
\frac{1}{n}
\int_{\tilde T_{i-1}}^{\tilde T_i}||\dot z(t)||^2\dt
\right]^\2  .  $$
If we compare the last term with the kinetic energy of $z$, we get that